\pdfoutput=1
\documentclass{amsart}

\usepackage{amsmath, amssymb, amsthm, fullpage, tikz, tikz-cd, pgf, enumerate, float, subcaption, bm}
\usepackage{pgfplots}
\usepackage{mathrsfs, mathabx}
\usetikzlibrary{arrows}
\usepackage[T1]{fontenc}
\usetikzlibrary{decorations.markings,patterns,arrows.meta}

\usepackage{newclude}
\usepackage{mathrsfs}
\usepackage{graphicx}
\usepackage{float}
\usepackage[english]{babel}
\usepackage[utf8]{inputenc}
\usepackage{csquotes}
\usepackage{fancyhdr}
\usepackage{verbatim}
\usepackage{changepage}
\usepackage{bbm}
\usepackage[miktex]{gnuplottex} 
\usepackage[margin=1.25in]{geometry}
\usepackage{relsize}
\usepackage{hyperref}

\theoremstyle{definition}
\newtheorem{definition}{Definition}[section]

\newtheorem{remark}[definition]{Remark}

\theoremstyle{plain}
\newtheorem{theorem}{Theorem}
\newtheorem{lemma}[definition]{Lemma}
\newtheorem{prop}[definition]{Proposition}
\newtheorem{cor}[definition]{Corollary}

\theoremstyle{remark}

\numberwithin{equation}{section}

\newcommand{\eq}[1]{\begin{align*} #1 \end{align*}}
\newcommand{\eqq}[1]{\begin{align} #1 \end{align}}
\newcommand{\scr}[1]{\mathscr{#1}}
\newcommand{\mc}[1]{\mathcal{#1}}

\newcommand{\R}{\mathbb{R}}
\newcommand{\C}{\mathbb{C}}

\newcommand{\Q}{\mathbb{Q}}
\newcommand{\Z}{\mathbb{Z}}

\newcommand{\set}[1]{\left\{#1\right\}}

\newcommand{\vspan}{\operatorname{span}}
\newcommand{\bb}[1]{\mathbb{#1}}

\newcommand{\ev}{\operatorname{ev}}
\newcommand{\ind}{\operatorname{ind}}
\newcommand{\CM}{\mathit{CM}}
\newcommand{\CS}{\mathit{CS}}
\newcommand{\HM}{\mathit{HM}}
\newcommand{\HS}{\mathit{HS}}
\newcommand{\pin}{\operatorname{Pin}}
\newcommand{\SWF}{\operatorname{SWF}}
\newcommand{\spinc}{spin$^c$ }
\newcommand{\Fps}{R}
\newcommand{\fps}{r}
\newcommand{\e}{\varepsilon}
\newcommand{\hess}{\operatorname{Hess}}
\newcommand{\gc}{\text{gC}}
\newcommand{\elc}{\text{elC}}
\newcommand{\agc}{\text{agC}}
\newcommand{\loc}{\text{loc}}
\newcommand{\inner}[1]{\left\langle#1\right\rangle}
\newcommand{\im}{\operatorname{im}}

\newcommand{\gr}{\operatorname{gr}}
\newcommand{\xgcq}{\mathcal X^{\gc}_{\frak q}}
\newcommand{\xgcql}{\mathcal X^{\gc}_{\frak q^\lambda}}

\newcommand{\xgcsq}{\mathcal X^{\gc,\sigma}_{\frak q}}
\newcommand{\xagcsq}{\mathcal X^{\agc,\sigma}_{\frak q}}
\newcommand{\xgcsql}{\mathcal X^{\gc,\sigma}_{\frak q^\lambda}}
\newcommand{\xagcsql}{\mathcal X^{\agc,\sigma}_{\frak q^\lambda}}

\begin{document}

\title{The equivalence of two Pin(2)-equivariant Seiberg-Witten Floer homologies}

\author{Nikhil Pandit}
\address{Department of Mathematics, Stanford University}
\email{npandit0@stanford.edu}

\begin{abstract}
We show that for a rational homology 3-sphere $Y$ equipped with a self-conjugate spin$^c$-structure $\mathfrak s$, the $\pin(2)$-equivariant monopole Floer homology of $(Y,\mathfrak s)$, as defined by Lin, is isomorphic to the $\pin(2)$-equivariant Seiberg-Witten Floer homology of $(Y,\mathfrak s)$ defined by Manolescu. 
\end{abstract}

\maketitle

\section{Introduction}\label{sec:intro}
Let $Y$ be a 3-manifold with a \spinc structure $\mathfrak s$.
In this paper, we are interested in two different frameworks for constructing Seiberg-Witten Floer homology of the pair $(Y,\mathfrak s)$. Firstly, there is Kronheimer and Mrowka's construction of monopole Floer homology in \cite{m3m}, based around achieving transversality through \textit{perturbations} of the Seiberg-Witten equations. This approach is very general: it is defined for all 3-manifolds $Y$. Monopole Floer homology comes in three flavors $\widecheck{\HM}_\ast(Y,\mathfrak s)$, $\widehat{\HM}_\ast(Y,\mathfrak s)$, and $\overline{\HM}_\ast(Y,\mathfrak s)$, which are modules over the ring $\Z[U]$ (where $\deg U=-2$) and are related by an exact sequence.

On the other hand, when $Y$ is a rational homology sphere, Manolescu's approach in \cite{man_orig}, based around finite-dimensional approximation, produces the Seiberg-Witten Floer spectrum $\SWF(Y,\mathfrak s)$, an $S^1$-equivariant suspension spectrum. Although this invariant is less general than monopole Floer homology, as it is not defined for all $3$-manifolds $Y$, it has a number of advantages when it \textit{is} defined: it can sometimes carry more information than the monopole Floer homology alone, and it interacts more simply with covers (see e.g. \cite{covering}). 
Accordingly, there have been several extensions of Manolescu's original construction to 3-manifolds $Y$ with $b_1(Y)>0$. When $b_1(Y)=1$ and the spin$^c$-structure is non-torsion, Kronheimer and Manolescu \cite{kronheimerPeriodicFloerProspectra2014} produce an invariant in the form of a periodic pro-spectrum. For 3-manifolds $Y$ with $b_1(Y)>0$,  Khandawit-Lin-Sasahira \cite{khandhawitUnfoldedSeibergWitten2018}  produce an ``unfolded'' Seiberg-Witten Floer spectrum. Finally, for 3-manifolds $Y$ for which the triple cup product on $H^1(Y;\Z)$ vanishes, Sasahira-Stoffregen \cite{sasahiraSeibergWittenFloerSpectra2021} produce a generalization of Manolescu's invariant for torsion spin$^c$-structures $\mathfrak s$.

In \cite{lm}, Lidman and Manolescu unify the monopole Floer and Seiberg-Witten Floer frameworks by proving an isomorphism of $\Z[U]$-modules between monopole Floer homology and the (reduced) $S^1$-equivariant homology of the spectrum, \eqq{
    \widecheck{\HM}_\ast(Y,\mathfrak s)\simeq \widetilde{H}^{S^1}_\ast(\text{SWF}(Y,\mathfrak s)).\label{eq:lm-iso}
    }

Note that the equivariant homology on the right-hand side is naturally a module over the ring $H^\ast(BS^1;\Z)=\Z[U]$.

Now suppose that $\mathfrak s$ is a self-conjugate \spinc structure on $Y$. In this case, there are equivariant refinements of both frameworks which incorporate the additional $\Z/2\Z$ symmetry. On the spectrum side, the extra $\Z/2\Z$ symmetry automatically produces a $\pin(2)$-action on the spectrum $\text{SWF}(Y,\mathfrak s)$, which was used by Manolescu \cite{triang} to produce a disproof of the triangulation conjecture in dimensions 5 and greater. On the monopole Floer homology side, Lin \cite{lin} constructs a refinement $\widecheck{\HS}_\ast(Y,\mathfrak s)$ (and similarly for the other flavors) by reworking Kronheimer and Mrowka's construction in the Morse-Bott setting. The invariant $\widecheck{\HS}_\ast(Y,\mathfrak s)$ is a module over the ring $\mc R=\bb F[Q,V]/(Q^3)$, where $\bb F=\Z/2\Z$, $Q$ has degree $-1$, and $V$ has degree $-4$. This invariant can be substituted for Manolescu's invariant to produce a corresponding disproof of the triangulation conjecture.

The purpose of this paper is to extend the isomorphism (\ref{eq:lm-iso}) to these $\pin(2)$-equivariant refinements, unifying the two corresponding approaches to the triangulation conjecture:

\begin{theorem}
Let $Y$ be a rational homology sphere, and let $\mathfrak s$ be a self-conjugate spin$^c$ structure on $Y$. Then there is an isomorphism of $\mc R$-modules
\eq{
    \widecheck{\HS}_\ast(Y,\mathfrak s)\simeq \widetilde{H}^{\pin(2)}_\ast(\operatorname{SWF}(Y,\mathfrak s);\bb F).}
Here $\widetilde{H}_\ast^{\pin(2)}$ is reduced $\pin(2)$-equivariant Borel homology, and the right-hand side is naturally a module over $H^\ast(B\pin(2);\bb F)=\mc R$.
\label{thm:main-thm}
\end{theorem}


The proof of Theorem \ref{thm:main-thm} proceeds similarly to the proof of the isomorphism (\ref{eq:lm-iso}) in \cite{lm}, and consists of the following steps.

First: 
\begin{enumerate}[(1)]
\item Show that $\widecheck{\HS}_\ast(Y,\mathfrak s)$ can be computed from configurations in global Coulomb gauge;
\item Show that finite-dimensional approximation can also be done in the presence of an equivariant perturbation, yielding a spectrum $\SWF_{\mathfrak q}(Y,\mathfrak s)$ which is $\pin(2)$-equivariantly stably homotopy equivalent to $\SWF(Y,\mathfrak s)$.
\end{enumerate}
These two steps are performed in \cite{lm}, and transfer over to the $\pin(2)$-equivariant setting essentially without change.

Next:
\begin{enumerate}
\item[(3)] Define an interpolating vector field $\mc X_\lambda$ on the global Coulomb slice $W$ by taking finite-dimensional projections of the nonlinear part of the perturbed Seiberg-Witten vector field $\mc X=\mc X_\infty$;
\item[(4)] Show that for $\lambda\gg0$, the stationary submanifolds of $\mc X_\lambda$ are in one-to-one correspondence with those of $\mc X$, and that the correspondence preserves gradings (after a suitable grading shift).
\item[(5)] Show that for $\lambda\gg 0$, the interpolating vector field $\mc X_\lambda$ defines a chain complex $\widecheck{\CS}^\lambda_\ast(Y,\frak s)$.
\end{enumerate}

These three steps are broadly similar to the approach taken in \cite{lm}, but many particulars are different. The correspondence in step (4) is obtained using the implicit function theorem, as in \cite{lm}, but because the stationary points are now non-isolated, the implicit function theorem has to be applied in a different way. Similarly, part of step (5) involves proving that the Hessian at reducible stationary points is \textit{normally} hyperbolic rather than fully hyperbolic, as in the non-degenerate case; and the chain complex is defined using Morse-Bott homology rather than Morse homology. 

The final ingredient of the proof is the following:
\begin{enumerate}
\item[(6)] Show that there is an $\mc R$-linear isomorphism $\widecheck{\HS}^\lambda_\ast(Y,\frak s)\simeq \widecheck{\HS}_\ast(Y,\frak s)$ (in a sufficiently large grading range).
\end{enumerate}

This section of the proof differs substantially from the corresponding step in \cite{lm}. There, since stationary points and trajectories were isolated, it was possible to apply the inverse function theorem to produce an explicit correspondence between trajectories and approximate trajectories, showing that the counts involved in the differentials were equal. The inverse function theorem also served to show that the $\mc R$-module structure was preserved. 

In the current setting, however, trajectories are non-isolated and the differentials are no longer simple counts: this setup is not amenable to the implicit function theorem. Instead, we produce the isomorphism in step (6) by considering moduli spaces for a time-dependent flow interpolating between the flow of $\mc X_\lambda$ and the flow of $\mc X$. These moduli spaces define a chain homotopy equivalence between the approximate complex $\widecheck{\CS}^\lambda(Y,\frak s)$ and the true complex $\widecheck{\CS}(Y,\frak s)$. Similarly, we show that the actions of $Q\in\mc R$ on the approximate and true complexes are related by a chain homotopy equivalence, and similarly for the actions of $V\in\mc R$. Thus, where \cite{lm} produced an isomorphism between chain complexes, we produce only a quasi-isomorphism.

The paper is organized as follows: In \S\ref{finite-dimensional}, we discuss the finite-dimensional model of Morse-Bott homology we will use. In \S\ref{sec:SW-prelim}, we briefly set up the objects and notation relevant to the analysis in the following sections. \S\ref{sec:gradings} sets up and compares various notions of gradings for stationary points of $\mc X$ and $\mc X_\lambda$. Finally, \S\ref{stationary-points} tackles steps (4) and (5) above, \S\ref{sec:inter-mod-sp} proves step (6), and \S\ref{sec:main-proof} puts the pieces together to prove Theorem \ref{thm:main-thm}.

\subsection*{Acknowledgement}
I am deeply grateful to my advisor Ciprian Manolescu for suggesting this problem and for his substantial help and advice throughout the research process.
\section{Morse-Bott homology for manifolds with boundary}\label{finite-dimensional}
This section defines a model of Morse-Bott homology in finite dimensions, later specializing to spaces with actions of the group $\pin(2)$.
\subsection{Overview}
Let $X$ be a finite-dimensional manifold with boundary, and let $v$ be a vector field on $X$ which is tangent to $\partial X$. The purpose of this section is to show that if $v$ satisfies a certain Morse-Bott-Smale quasi-gradient condition, then for every isolated invariant set $\mathscr S$ of the flow of $v$, we can define a Morse-Bott chain complex $\widecheck{C}_\ast(X,v)[\mathscr S]$ whose homology computes the reduced homology of the Conley index associated to $\mathscr S$ and the flow of $v$. The generators of the chain complex $\widecheck{C}_\ast(X,v)[\mathscr S]$ 
will be maps $\sigma:\Delta\to B$ from a stratified space $\Delta$ into a stationary submanifold $B$ of $v$, 
and the differential $\partial\sigma$ will be a sum $\sum_{B'} \Delta\times_{\sigma,\ev_-} \mc M(B,B')$ of fibered products of
moduli spaces of trajectories from $B$ to another stationary submanifold $B'$. Therefore, it is necessary to start by defining a class of stratified spaces containing the moduli spaces $\mathcal M(B,B')$ and showing that the class is closed under fiber products. This is done in \cite{lin}, but we reproduce the basic definitions and results below.
\subsection{$\delta$-chains}

Let $B, B'$ be two stationary submanifolds of the vector field $v$. Because $v$ is tangent to $\partial X$, the compactified moduli space $\breve M^+(B,B')$ of broken flow lines of $v$ is not naturally a smooth manifold with corners (see Section 16.5 in \cite{m3m}). However, when $v$ satisfies a transversality condition, the compactified moduli space will have the structure of an \textit{abstract $\delta$-chain}. We start with the following two definitions to prepare the definition of an abstract $\delta$-chain:
\begin{definition}[Definition 3.5.1 in \cite{lin}]
Let $(Q,q_0)$ be a pointed topological space. Let $\pi:S\to Q$ be a continuous map, and suppose $S_0\subset \pi^{-1}(q_0)$. We say that $\pi$ is a \textit{topological submersion along $S_0$} if for every $s_0\in S_0$ there is a neighborhood $U\subset S$ and a neighborhood $Q'\subset Q$ of $q_0$ admitting a homeomorphism \eq{
	(U\cap S_0)\times Q'\to U}
which commutes with $\pi$.
\label{top-submersion}
\end{definition}
\begin{definition}[Definition 3.6.13 in \cite{lin}]
Let $N$ be a $d$-dimensional space stratified by manifolds and $M^{d-1}$ a union of $(d-1)$-dimensional components of the stratification. We say that $N$ has a \textit{codimension-$c$ $\delta$-structure} along $M^{d-1}$ if there is an open set $W\supset M^{d-1}$ admitting a topological embedding $j: W\to EW$ and a map \eq{
	\mathbb S=(S_1,\dots, S_{c+1}):EW\to (0,\infty]^{c+1}}
	with the following properties:
\begin{enumerate}[(1)]
\item the fiber along $\infty$ is identified with $j(M^{d-1})$ and $\mathbb S$ is a topological submersion along it;
\item the subset $j(W)\subset EW$ is the zero set of a continuous map $\delta: EW\to \Pi^c$ (where $\Pi^c\subset \R^{c+1}$ denotes the hyperplane $\set{\delta\in \R^{c+1}\mid\sum\delta_i=0}$);
\item if $e\in EW$ has $S_{i_0}=\infty$ for some index $i_0$, then $\delta_{i_0}\le 0$ (with equality only if $\mathbb S(e)=\infty$);
\item on the subset of $EW$ where all the $S_i$ are finite, $\delta$ is smooth and transverse to zero.
\end{enumerate}
\label{codim-c-structure}
\end{definition}
Now we can define an abstract $\delta$-chain: 
\begin{definition}[Definition 4.1.1 in \cite{lin}]
A topological space $N^d$ is a \textit{$d$-dimensional abstract $\delta$-chain} if it is a stratified space\eq{
	N^d\supset N^{d-1}\supset\dots\supset N^0\supset N^{-1}=\varnothing} of dimension $d$ with the following additional structures. We are given a finite partition of each stratum \eq{
	N^e\setminus N^{e-1}=\coprod_{i=1}^{m_e} M_i^e} for each $e=0,\dots, d$, and we call the closure of each $M_i^e$ an $e$-dimensional \textit{face}. We denote the top stratum of a face $\Delta$ by $\mathring\Delta$. The set of faces satisfies the following combinatorial property: whenever a codimension $e$ face $\Delta''$ is contained in a codimension $e-2$ face $\Delta$, there are exactly two codimension $e-1$ faces contained in $\Delta$ that contain $\Delta''$.
	
Furthermore, each pair of faces $\Delta'\subset \Delta$ has an associated finite set $N(\Delta, \Delta')$, together with a subset $O(\Delta, \Delta')\subset N(\Delta,\Delta')$, where the $O$ is for obstructed, and their local structure satisfies the following properties. There is an open neighborhood $\breve W(\Delta,\Delta')$ of $\mathring \Delta'$ inside $\Delta$ together with a topological embedding $j$ inside a space $E\breve W(\Delta,\Delta')$ endowed with a topological submersion \eq{
	\mathbf{S}:E\breve W(\Delta,\Delta')\to(0,\infty]^{N(\Delta,\Delta')}}
called the \textit{local thickening}, such that:
\begin{enumerate}[(1)]
\item The map $\mathbf{S}$ is a topological submersion along the fiber over $\mathbf{\infty}$, and this is identified with $\mathring \Delta'$;
\item The image $j(\breve W)\subset E\breve W$ is the zero set of a map \eq{
	\delta: E\breve W\to\R^{O(\Delta,\Delta')}}
	vanishing along the fiber of $\mathbf{S}$ over $\mathbf{\infty}$;
\item Calling $\breve W^0\subset \breve  W$ and $E\breve W^0\subset E\breve W$ the subsets where none of the components of $\mathbf{S}$ is infinite, the restriction of $j$ to $\breve W^0$ is a smooth embedding, and the restriction of $\delta$ to $E\breve W^0$ is transverse to zero.
\end{enumerate}

This collection of local thickenings is \textit{compatible} in the following sense: Whenever we have three faces $\Delta''\subset\Delta'\subset \Delta$, we have a canonical inclusion \eq{
	N(\Delta',\Delta'')\hookrightarrow N(\Delta,\Delta'')} 
	with identifications \eq{
	O(\Delta',\Delta'')=O(\Delta,\Delta'')\cap N(\Delta',\Delta'')}
and \eq{
	(0,\infty]^{N(\Delta',\Delta'')}=\{(x_1,\dots,x_{N(\Delta,\Delta'')})\mid x_\alpha=\infty \text{ if } \alpha\notin N(\Delta',\Delta'')\}\subset (0,\infty]^{N(\Delta,\Delta'')}}
so that we also have the identifications \eq{
	E\breve W(\Delta',\Delta'')&\equiv E\breve W(\Delta,\Delta'')\cap\mathbf S^{-1}\left((0,\infty]^{N(\Delta',\Delta'')}\times\{\infty\}\right)\\
	\mathbf S(\Delta',\Delta'')&\equiv \mathbf S(\Delta,\Delta'')|_{E\breve W(\Delta',\Delta'')}:E\breve W(\Delta',\Delta'')\to(0,\infty]^{N(\Delta',\Delta'')}\\
	\delta(\Delta',\Delta'')&\equiv \delta(\Delta,\Delta'')|_{E\breve W(\Delta',\Delta'')}: E\breve W(\Delta',\Delta'')\to \R^{O(\Delta',\Delta'')}.}
	(In the last line, we are implicitly stating that the restriction of $\delta(\Delta,\Delta'')$ has image contained in the subspace corresponding to $\R^{O(\Delta',\Delta'')}$ under the identifications above.) Finally, whenever $\Delta'\subset \Delta$ has codimension one, for each face $\Delta$ the data above defines a codimension $c$ $\delta$-structure along $\Delta'$ in the sense of Definition \ref{codim-c-structure}, in the sense that there is an identification of $\R^{O(\Delta,\Delta')}$ with the subspace $\Pi^c\subset \R^{c+1}$.
\label{abstract-delta-chain}
\end{definition}

\begin{definition}[Definition 4.1.6 in \cite{lin}]
Let $X$ be a smooth manifold without boundary. A \textit{$\delta$-chain} in $X$ is a pair $\sigma=(\Delta, f)$ where 
\begin{enumerate}[(1)]
\item $\Delta$ is an abstract $\delta$-chain and $f:\Delta\to X$ is a continuous map;
\item the restriction of $f$ to each stratum of $\Delta$ is a smooth map;
\item for each pair of faces $\Delta'\subset \Delta''$, the map extends to a continuous map $Ef(\Delta',\Delta'')$ on the local thickening $E\breve W(\Delta',\Delta'')$ which is smooth in a neighborhood of $\breve W(\Delta,\Delta'')$;
\item the collection of extensions to the local thickenings is compatible, in the sense that for every triple $\Delta\supset\Delta'\supset\Delta''$ the map $Ef(\Delta',\Delta'')$ is the restriction of $Ef(\Delta,\Delta'')$ under the identification of Definition \ref{abstract-delta-chain}.
\end{enumerate}
\label{delta-chain}
\end{definition}

\begin{lemma}[Lemma 4.1.8 in \cite{lin}]
The fiber product of two \textit{transverse} $\delta$-chains in $X$ has the natural structure of a $\delta$-chain in $X$.
\end{lemma}

\begin{definition}[p. 96 in \cite{lin}]
For any countable family $\mathcal F$ of $\delta$-chains in $X$, we write $C^{\mathcal F}_d(X)$ for the free vector space over $\bb F=\Z/2\Z$ generated by $d$-dimensional $\delta$-chains in $X$ which are transverse to the family $\mathcal F$, modulo isomorphism and modulo ``negligible chains'' (Definition 4.1.11 in \cite{lin}). There is a boundary map $\partial_0:C_d^{\mathcal F}\to C_{d-1}^{\mathcal F}$, defined on a $\delta$-chain $(\Delta,\sigma)$ by summing the restrictions of $\sigma$ to the codimension-1 faces of $\Delta$.
\label{crit-chains}
\end{definition}

\begin{remark}
As in Remark 4.1.12 in \cite{lin}, $C_d^{\mathcal F}(X)$ is trivial for $d\ge \operatorname{dim}(X)+2$.
\label{rmk:negl-chains}
\end{remark}
\begin{prop}[Proposition 4.1.13 in \cite{lin}]
The homology of the complex $(C_\ast^{\mathcal F}(X),\partial_0)$ is canonically isomorphic to the singular homology $H_\ast(X; \bb F)$.
\label{lip-equivalence}
\end{prop}

\subsection{Morse-Bott quasi-gradients}
Although there are standard definitions of Morse-Bott homology for Morse-Bott gradient fields $v=\text{grad } f$, we will be forced to deal with vector fields which are \textit{approximations} to the Seiberg-Witten gradient flow, and hence are not themselves gradients. However, these approximate vector fields will still satisfy the following ``quasi-gradient'' condition, which extends Definition 2.3.2 in \cite{lm} to the Morse-Bott setting.
\begin{definition}
Let $X$ be a smooth manifold, possibly with boundary. A smooth vector field $v$ on $X$ is called a Morse-Bott quasi-gradient if
\begin{enumerate}[(i)]
\item $v$ is tangent to $\partial X$ along $\partial X$.
\item The stationary points of $v$ form closed smooth manifolds, each of which is contained entirely in $\partial X$ or in $X-\partial X$.
\item For each stationary submanifold $B$, the derivative $dv$ is \textit{normally hyperbolic}; i.e. for each $x\in B$, the complexification of $(dv)_x:T_xX/T_xB\to T_xX/T_xB$ has no eigenvalues with real part 0. \label{item-normhyp}
\item There exists a smooth function $f:X\to\R$ such that $df(v)\ge 0$ at all $x\in X$, with equality holding if and only if $x$ is a stationary point of $v$. \label{item-tamingfn}
\end{enumerate}
\label{morse-bott-quasigradient-def}
\end{definition}

A stationary submanifold $B$ of a Morse-Bott quasi-gradient has an index $\ind(B)$, namely the number of eigenvalues of $dv$ on $B$ with negative real part. Our task is to define a version of Morse-Bott homology for Morse-Bott quasi-gradients.
\par 
The following lemma, which follows from Theorem 2 in \cite{pughshub}, shows that condition (\ref{item-normhyp}) above guarantees that a Morse quasi-gradient is conjugate to a linear flow near a stationary submanifold:
\begin{lemma}
Assume that a smooth vector field $v$ on a Riemannian manifold $X$ has a compact stationary submanifold $B$ on which $dv$ is normally hyperbolic, as in item (\ref{item-normhyp}) above. Then there is a homeomorphism from a neighborhood of $B$ to a neighborhood of the zero section in the normal bundle $NB$ which conjugates the flow of $v$ to the flow of the linear part of $v$.
\label{morsebott-locnorm-lemma}
\end{lemma}

On the other hand, condition (\ref{item-tamingfn}) in Definition \ref{morse-bott-quasigradient-def} guarantees that flow lines have limits in our setting (Lemma \ref{limits-exist-noneq} below). This lemma is a crucial ingredient in showing that the stable and unstable manifolds of $v$ are embedded (Proposition \ref{prop:stable-mfld-embedded} below).

\begin{lemma}
Let $v$ be a Morse-Bott quasi-gradient vector field on a compact, oriented, smooth manifold $X$. Let $\gamma:\R\to X$ be a flow line of $v$ (i.e. $\frac{d\gamma}{dt}+v(\gamma)=0$). Then $\lim_{t\to\infty}\gamma(t)$ and $\lim_{t\to -\infty}\gamma(t)$ exist, and they are both stationary points of $v$.
\label{limits-exist-noneq}
\end{lemma}
\begin{proof}
Let $\mathfrak C=\{B_1,\dots, B_k\}$ be the set of connected critical submanifolds of $v$. By Lemma \ref{morsebott-locnorm-lemma}, for each $B\in\mathfrak C$ we may choose a neighborhood $\mc N_B$ and a homeomorphism $\phi_B:\mc N_B\to D^k\times D^{n-k}\times B$ conjugating $v$ with the vector field \eq{
    \tilde v(x_1,\dots, x_n,b)=(-x_1,\dots,-x_k,x_{k+1},\dots, x_n,0).}
We can assume the neighborhoods $\mc N_B$ are disjoint. 

\begin{figure}[h!]
\centering
\begin{tikzpicture}[scale=0.8, decoration={markings, 
    mark=at position 0.2 with {\arrow{>}},
    mark=at position 0.5 with {\arrow{>}},
    mark=at position 0.8 with {\arrow{>}}
    }]

\draw[thick,gray,->] (0,-4) -- (0,4) node[above right] {};
\draw[thick,gray,->] (-4,0) -- (4,0) node[right] {};

\draw[thick] (-4,-4) rectangle (4,4);

\draw[very thick,teal] (0.5,4) -- (-0.5,4) node[above right] {$A_{\text{in}}$};
\draw[very thick,teal] (0.5,-4) -- (-0.5,-4) node[below right] {$A_{\text{in}}$};
\draw[very thick,red] (4,-0.5) -- (4,0.5) node[below right] {$A_{\text{out}}$};
\draw[very thick,red] (-4,-0.5) -- (-4,0.5) node[below left] {$A_{\text{out}}$};

\draw[thick,blue,postaction={decorate}] plot[smooth, domain=0.5:4, samples=10] (\x,{2/\x});
\draw[thick,blue,postaction={decorate}] plot[smooth, domain=0.5:4, samples=10] (\x,{-2/\x});
\draw[thick,blue,postaction={decorate}] plot[smooth, domain=-0.5:-4, samples=10] (\x,{2/\x});
\draw[thick,blue,postaction={decorate}] plot[smooth, domain=-0.5:-4, samples=10] (\x,{-2/\x});

\draw[blue,postaction={decorate}] plot[smooth, domain=0.25:4, samples=10] (\x,{1/\x});
\draw[blue,postaction={decorate}] plot[smooth, domain=0.25:4, samples=10] (\x,{-1/\x});
\draw[blue,postaction={decorate}] plot[smooth, domain=-0.25:-4, samples=10] (\x,{1/\x});
\draw[blue,postaction={decorate}] plot[smooth, domain=-0.25:-4, samples=10] (\x,{-1/\x});

\fill[yellow, opacity=0.2] 
(0.5, 4) -- (-0.5, 4) -- 
plot[smooth, domain=-0.5:-4, samples=14] (\x, {2/\x}) --
plot[smooth, domain=-4:-0.5, samples=14] (\x, {-2/\x}) --
(-0.5, -4) -- (0.5, -4) -- 
plot[smooth, domain=0.5:4, samples=14] (\x, {-2/\x}) --
plot[smooth, domain=4:0.5, samples=14] (\x, {2/\x}) -- 
cycle;

\node[right,blue] at (1,2) {$A_{\text{flow}}$};
\node[above,black] at (1.5,-0.05) {$D^k$};
\node[right,black] at (-0.1,0.7) {$D^{n-k}$};

\end{tikzpicture}
\caption{Neighborhood $A$ of $0$ (shaded) in $D^k\times D^{n-k}$}
\label{fig:morse-nbhd}

\end{figure}
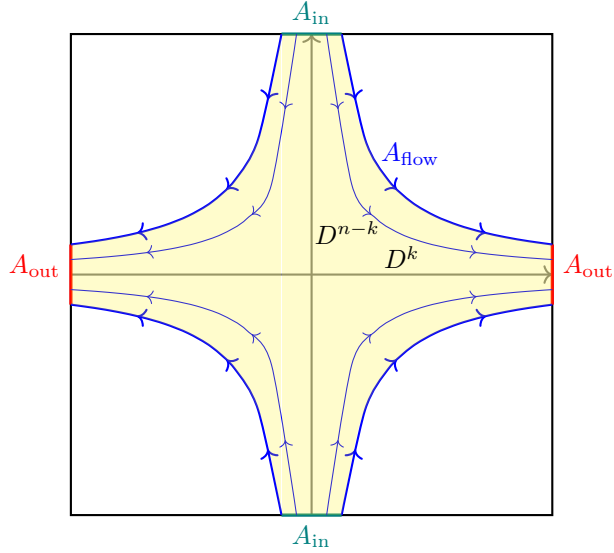

For sufficiently small $\eta>0$, consider the subset $A\subset D^k\times D^{n-k}$ consisting of all points on flow lines of $-\tilde v$ starting on $D^k(\eta)\times S^{n-k}(1)\times B$ or at $\{0\}\times\{0\}\times B$ (Figure \ref{fig:morse-nbhd}). The boundary of $A$ consists of the three pieces $A_{\text{in}}=D^k(\eta)\times S^{n-k}(1)\times B$, $A_{\text{out}}=S^k(1)\times D^{n-k}(\eta)\times B$, and $A_{\text{flow}}$ consisting of points on flow lines between $\partial A_{\text{in}}$ and $\partial A_{\text{out}}$.

Let $\mc A_B$, $\mc A_{B,\text{in}}$, $\mc A_{B,\text{out}}$ be the preimages of $A$, $A_{\text{in}}$, and $A_{\text{out}}$ under $\phi_B$. Let $c_B$ be the value of $f$ on $B$. Since the subset $A_{\text{core}}=\{0\}\times S^{n-k}(1)\times B\subset A_{\text{in}}$ flows into $B$ under the flow of $-\tilde v$, and since $f$ decreases on trajectories of $-v$, the value of $f\circ \phi_B^{-1}$ on any point of $A_{\text{core}}$ is greater than $c_B$. By compactness of $S^{n-k}(1)\times B$, we can choose $\eta>0$ small enough that the value of $f$ on any point of $\mc A_{B,\text{in}}$ is greater than $c_B$. Similarly, we can ensure that the value of $f$ on $\mc A_{B,\text{out}}$ is less than $c_B$.

Let $\mc N_B$ be the interior of $\mc A_B$. By construction, flow lines must enter $\mc N_B$ through $\mc A_{B,\text{in}}$ and exit through $\mc A_{B,\text{out}}$. Furthermore, $f$ decreases along flow lines, and $f|_{\mc A_{B,\text{in}}}>c_B>f|_{\mc A_{B,\text{out}}}$. It follows that if a flow line enters $\mc N_B$, it either remains in $\mc N_B$ for all future time, or it exits and never returns.

Let $\mathcal N=\bigcup_{B\in\mathfrak C} \mc N_B$, and choose $\e>0$ so that $df(v)\ge \e$ on $X-\mc N$. Let $\gamma:\R\to X$ be a flow line of $v$, i.e. a solution to the equation $d\gamma/dt + v(\gamma(t))=0$. If $\gamma$ stayed outside $\mc N$ for all $t\ge t_0$, we would have \eq{
    f(\gamma(T))= f(\gamma(t_0))-\int_{t_0}^T df(v(\gamma(t)))\,dt\le f(\gamma(t_0))-\e (T-t_0)} 
for all $T\ge t_0$, which is impossible since $f$ is bounded. Therefore $\gamma$ either remains in $\mc N$ or leaves and returns to $\mc N$ infinitely many times. However, since there are finitely many $\mc N_B$'s and $\gamma$ can never return to a neighborhood $\mc N_B$ once it has left, $\gamma$ must remain in $\mc N$, and hence some $\mc N_B$, for all $t$ sufficiently large. It follows that $\lim_{t\to\infty} \gamma(t)$ exists, as flow lines of $\tilde v$ which remain in $A$ for all positive time have limits at $\infty$. We can argue similarly to conclude that $\lim_{t\to -\infty} \gamma(t)$ exists as well.
\end{proof}
\begin{remark}
The proof of Lemma \ref{limits-exist-noneq} is similar to the proof of Lemma 2.8 in \cite{lm}, which was missing the argument involving controlling the values of $f$ on $\mc A_{B,\text{in}}$ and $\mc A_{B,\text{out}}$.
\end{remark}

Given $x\in X$, let $\gamma_x$ be the flow line of the Morse-Bott quasi-gradient $v$ satisfying $\gamma_x(0)=x$. For a stationary submanifold $B\subset X$ of $v$, we define \eq{
    W^s(B)&=\set{x\in X \,\Big| \lim_{t\to\infty} \gamma_x(t)\in B}\\
    W^u(B)&=\set{x\in X \,\Big| \lim_{t\to-\infty} \gamma_x(t)\in B}.}
We can also define $W^s(p)$, $W^u(p)$ in the same way when $p$ is a (possibly non-isolated) stationary point of $v$.
\begin{prop}
The stable and unstable sets $W^s(B), W^u(B)$ above are embedded submanifolds of $X$.
\label{prop:stable-mfld-embedded}
\end{prop}
\begin{proof}
According to Theorem A.9 in \cite{ab}, an open neighborhood of $B$ in $W^s(B)$ can be identified with a neighborhood of the zero section in the stable normal bundle $\nu^-\to B$, which consists of the generalized eigenspaces of $dv$ corresponding to eigenvalues with negative real part. With this local description in hand, we can construct a homeomorphism $E^s:\nu^-\to W^s(B)$ exactly as in the Morse case (see Theorem 4.20 in \cite{bhmorse}). The key point in that theorem is the fact that flow lines have limits at $\infty$, which is proved in the current setting as Lemma \ref{limits-exist-noneq}. The case of $W^u(B)$ is analogous.
\end{proof}

\subsection{The Morse-Bott-Smale condition}
In this section, we define moduli spaces of trajectories of a Morse-Bott quasi-gradient under an additional transversality condition.
\begin{definition}
We say that a Morse-Bott quasi-gradient $v$ is \textit{Morse-Bott-Smale} if for all stationary submanifolds $B_-,B_+\subset X$ of $v$, the intersection $W^u(B_-)\pitchfork W^s(B_+)$ is transverse in $X$. (If $X$ has boundary and both $W^u(B_-)$ and $W^s(B_+)$ are contained in $\partial X$, we instead require $W^u(B_-)$ and $W^s(B_+)$ to intersect transversely in $\partial X$.)
\label{basic-morse-bott-smale-def}
\end{definition}
When $v$ is Morse-Bott-Smale, we can define $\breve M(B_-,B_+)=(W^u(B_-)\cap W^s(B_+))/\R$, the moduli space of unparametrized flow lines of $v$ from $B_-$ to $B_+$. There are also endpoint evaluation maps $\ev_-:\breve M(B_,B+)\to B_-$, $\ev_+:\breve M(B_-,B_+)\to B_+$ which take a trajectory to its asymptotic start/end point.

\begin{remark}
Definition \ref{basic-morse-bott-smale-def} differs from the Smale condition used in  \cite{ab} (Section 3.2), where one imposes the stronger condition that for every $p\in B_-$, the intersection $W^u(p)\pitchfork W^s(B_+)$ is transverse. There, this asymmetric definition is necessary to ensure that the map $\ev_-$ defined above is a locally trivial fiber bundle. However, in the $\delta$-chain formulation of singular homology, it is no longer necessary to impose this condition. 
\end{remark}
In the remainder of this section, we express the Morse-Bott-Smale transversality condition equivalently as surjectivity of a certain linear operator.

We begin by defining the weighted Sobolev spaces necessary for the discussion. Fix two stationary submanifolds $B_-$ and $B_+$, and set \eq{
\mc P_k(B_-,B_+)=\set{\gamma\in L^2_k(\R, X)\,\Big|\, \lim_{t\to\pm\infty}\gamma(t)\in B_{\pm}}.} 
Choose $\delta>0$ small enough that for all $p\in B_-\cup B_+$ and for all nonzero eigenvalues $\lambda$ of $dv_p$, $|\text{Re }\lambda|>\delta$. Choose once and for all a function $f:\R\to\R$ with $f(t)=e^{\delta|t|}$ for $|t|\gg 0$, and for any $L^2_j$ vector field $w$ along a path $\gamma\in P_k(B_-,B_+)$, define $\|w\|_{j,\delta}=\|fw\|_j$. Finally, we set $T_{\gamma,j,\delta} \mc P(B_-,B_+)$ to be the completion of the tangent space to $\mc P_k(B_-,B_+)$ at $\gamma$ with respect to the $L^2_{j,\delta}$ norm $\|\cdot\|_{j,\delta}$. (Note that multiplication by $f$ defines an isometry $M_f:L^2_{k,\delta}\to L^2_k$.)

Linearizing the flow operator $F:\mc P_k(B_-,B_+)\to T_{k-1} \mc P(B_-,B_+)$ \eq{
    F(\gamma)=\frac{d\gamma}{dt} + v(\gamma),}
    gives rise to a linearized operator $L_\gamma:T_{\gamma,j,\delta} \mc P(B_-,B_+)\to T_{\gamma,j-1,\delta} \mc P(B_-,B_+)$, \eqq{
    L_\gamma&=\frac{Dw}{dt}+Dv(w).\label{l-gamma-eq}}
Here $D$ is any connection on $TX$ defined in a neighborhood of the image of $\gamma$. The operator $L_\gamma$ is Fredholm (cf. Proposition 3.3.4 in \cite{lin}). As in Lemma 2.1 in \cite{lm}, the particular choice of the connection $D$ only changes the operator $L_\gamma$ by a compact operator, and when $\gamma$ is a flow line $L_\gamma$ is independent of the choice of $D$.

\begin{prop}
The vector field $v$ is Morse-Bott-Smale if and only if for any critical submanifolds $B_-$, $B_+$ of $v$ and for any flow line $\gamma$ from $B_-$ to $B_+$, the operator $L_\gamma$ from (\ref{l-gamma-eq}) is surjective. 
\end{prop}
\begin{proof}
The proof closely mirrors Proposition 10.2.8 in \cite{ad}. Fix $B_-$, $B_+$, and $\gamma$ as in the proposition, and let $b_{\pm}=\ev_{\pm}(\gamma)$. Given $s_0, s\in \R$, we let $\Phi_{(s_0,s_1)}: T_{\gamma(s_0)} X \to T_{\gamma(s_1)}X$ be the resolvent of the linear ODE, i.e. the unique solution $w$ to $L_\gamma(w)=0$ with $w(s_0)=Y$ satisfies $w(s_1)=\Phi_{(s_0,s_1)}(Y)$. Set\eq{
    E^u(s_0)&=\set{Y\in T_{\gamma(s_0)}X\,\Big|\, \lim_{s\to -\infty} f(s)\Phi_{(s_0,s)}Y=0}\\
    E^s(s_0)&=\set{Y\in T_{\gamma(s_0)}X\,\Big|\, \lim_{s\to +\infty} f(s)\Phi_{(s_0,s)}Y=0}\\
    E^o_{\pm}(s_0)&=\set{Y\in T_{\gamma(s_0)}X\,\Big|\, \Phi_{(s_0,s)}Y\text{ stays bounded away from } 0\text{ and }\infty \text{ as }s\to \pm\infty}.} Then $E^u(s_0)=T_{\gamma(s_0)} W^u(b_-)$, $E^s(s_0)=T_{\gamma(s_0)} W^s(b_+)$, and
$\ker L_u$ is identified with $E^u(s_0)\cap E^s(s_0)$. (The spaces $E^o_\pm(s_0)$, meanwhile, are identified with $T_{b_{\pm}}B_{\pm}$.)

To investigate the cokernel of $L_\gamma$, we investigate the kernel of its $L^2$-formal adjoint \eq{
    L_\gamma^\ast=-\frac{D}{dt}+(Dv)^\ast.} It turns out that $\Psi_{(s_0,s_1)}:=\Phi_{(s_1,s_0)}*$ is the resolvent of $L_u^\ast$, for differentiating the identity \eqq{
    \inner{\Psi_{(s_0,s)}Z,\Phi_{(s_0,s)}Y}=\inner{Z,Y}\label{smale-identity}}
with respect to $s$ and writing $Z(s)=\Psi_{(s_0,s)}(Z)$, $Y(s)=\Phi_{(s_0,s)}(Y)$ yields \eq{
    \inner{\frac{DZ(s)}{ds}, Y(s)}+\inner{Z(s),\frac{DY(s)}{ds}}&=0\\
    \inner{\frac{DZ(s)}{ds}-(Dv)^\ast Z(s),Y(s)}=0\\
    \inner{L_\gamma^\ast Z(s),Y(s)}=0,} and since this is true for any $Y$ we conclude $L_\gamma^\ast Z(s)=0$.
    
Investigating \ref{smale-identity} more closely, we see that if $Z\in \ker L_\gamma^\ast$, then $\inner{Z,Y}=0$ whenever $\Phi_{(s_0,s)}Y$ remains bounded as $s\to\pm\infty$. This shows that \eq{
    \ker L_\gamma^\ast&\cong(E^u(s_0)+E^o_-(s_0))^\perp\cap (E^s(s_0)+E^o_+(s_0))^\perp\\
    &\cong (E^u(s_0)+E^o_-(s_0)+E^s(s_0)+E^0_+(s_0))^\perp.}
Since $E^u(s_0)+E^o_-(s_0)=T_{\gamma(s_0)} W^u(B_-)$ and $E^s(s_0)+E^o_+(s_0)=T_{\gamma(s_0)} W^s(B_+)$, we see that $L_\gamma^\ast$ is injective if and only if $W^u(B_-)$ intersects $W^s(B_+)$ transversely, as desired.

\end{proof}

\subsection{The Morse-Bott chain complex}
Recall that if $X$ is a manifold with boundary, and if $f$ is a Morse function on $X$ whose gradient is parallel to $\partial X$, Kronheimer and Mrowka defined a Morse chain complex involving three types of critical points: \textit{interior}, \textit{boundary-stable}, and \textit{boundary-unstable} critical points. Interior critical points lie in $X-\partial X$, while critical points $p$ in the boundary are said to be boundary-stable or boundary-unstable if the normal vector to the boundary at $p$ is stable or unstable, respectively, for the flow of $-\text{grad }f$. A trajectory from a critical point $p$ to a critical point $q$ is said to be \textit{boundary-obstructed} if $p$ is boundary-stable and $q$ is boundary-unstable.

Let $v$ be a Morse-Bott-Smale quasi-gradient field on $X$, and fix an isolated invariant set $\mathscr S$ of the flow of $v$, in the sense of Definition 2.4.1 in \cite{lm}. Stationary submanifolds of $v$ now come in the same three types (interior, boundary-stable, boundary-obstructed). For two stationary submanifolds $B, B'\subset \mathscr S$ of $v$, let $\breve M[\mathscr S](B, B')$ be the moduli space of unparametrized flow lines of $v$ from $B$ to $B'$ which lie entirely within $\mathscr S$. Write $\ev_-:\breve M(B,B')\to B$, $\ev_+:\breve M(B,B')\to B'$ for the maps which take a trajectory to its asymptotic start/end point.

\begin{prop}
The moduli space $\breve M[\mathscr S](B,B')$ has a compactification \eq{
	\breve M^+[\mathscr S](B,B')=\bigcup \breve M[\mathscr S](B,B_1)\times_{B_1} \breve M[\mathscr S](B_1,B_2)\times_{B_2}\dots\times_{B_n} \breve M[\mathscr S](B_n,B'),} where we take the union over all sequences $(B_1,\dots, B_n)$ of critical submanifolds, and where the fiber products are taken with respect to the maps $\ev_{\pm}$.
	\label{compactification}
\end{prop}
\begin{proof}
This is stated underneath Theorem 4.9 in \cite{bh} for the case that $v$ is a gradient. However, it is also true when $v$ is a Morse-Bott-Smale quasi-gradient, since the underlying analysis in Appendix \S A.3 of \cite{ab} is done under the sole assumption that $v$ has hyperbolic stationary manifolds.
\end{proof}
The following proposition is analogous to Proposition 3.6.14 in \cite{lin}, and the dimension calculation is analogous to Proposition 3.3.12 from the same source.
\begin{prop}
The compactified moduli space $\breve M^+[\mathscr S](B,B')$ has the structure of an abstract $\delta$-chain. Its dimension is \eq{
	\operatorname{dim}\left(M^+[\mathscr S](B,B')\right)&= \begin{cases} \ind(B)-\ind(B')+\dim B-1 & (B,B')\text{ not boundary-obstructed}\\
	\ind(B)-\ind(B')+\dim B&(B,B')\text{ boundary-obstructed.}\end{cases}}
\label{prop:delta-chain-with-dim}
\end{prop}

Each compactified moduli space $\breve M^+[\mathscr S](B,B')$ determines a $\delta$-chain $(\breve M^+[\mathscr S](B,B'),\ev_-)$ in $X$ (and in fact in $B$). Let $\mathcal F=\set{(\breve M^+[\mathscr S](B,B'),\ev_-)\mid B,B'\in\text{Crit}(f)}$. Then we may define $C^{\mathcal F}_\ast(B)$ as in Definition \ref{crit-chains}.

Recall that if $X$ is a manifold with boundary, and if $f$ is a Morse function on $X$ whose gradient is parallel to $\partial X$, Kronheimer and Mrowka defined a Morse chain complex involving three types of critical points: \textit{interior}, \textit{boundary-stable}, and \textit{boundary-unstable} critical points. Interior critical points lie in $X-\partial X$, while critical points $p$ in the boundary are said to be boundary-stable or boundary-unstable if the normal vector to the boundary at $p$ is stable or unstable, respectively, for the flow of $-\text{grad }f$. A trajectory from a critical point $p$ to a critical point $q$ is said to be \textit{boundary-obstructed} if $p$ is boundary-stable and $q$ is boundary-unstable.

For our Morse-Bott-Smale vector field $v$, the stationary submanifolds of $v$ come in the same three types. Let $\mathfrak C^o[\mathscr S](v)$, $\mathfrak C^s[\mathscr S](v)$, and $\mathfrak C^u[\mathscr S](v)$ denote the collection of interior, boundary-stable, and boundary-unstable stationary submanifolds, respectively. For $\theta\in\set{o,s,u}$, define \eqq{
	C^\theta_k[\mathscr S] =\bigoplus_{i+j=k}\bigoplus_{\substack{B\in\mathfrak C^\theta[\mathscr S](v)\\ \ind(B)=j}}C^{\mathcal F}_i[\mathscr S](B).\label{breakup}}
	
Write $C_k[\mathscr S]=C^o_k[\mathscr S]\oplus C^s_k[\mathscr S]\oplus C^u_k[\mathscr S]$. There is a differential $\partial_{>0}: C_k[\mathscr S]\to C_{k-1}[\mathscr S]$ which takes a $\delta$-chain $\sigma:\Delta\to B$ to the sum $\sum_i \sigma\times_B \breve M^+[\mathscr S](B,B_i)$, where the summand is viewed as a $\delta$-chain in $B_i$. (Note that $\partial_{>0}$ splits as \eq{
	\partial_{>0}=\partial_1+\partial_2+\dots+\partial_N,} where $\partial_m$ takes chains in an index $j$ stationary submanifold to chains in an index $j-m$ stationary submanifold, and $N=\dim X$.) We set $\partial=\partial_0+\partial_{>0}$, where $\partial_0$ is as defined in Definition \ref{crit-chains}.
	For $\theta,\theta'\in \set{o,s,u}$, let $\partial^\theta_{\theta'}$ be the component of $\partial$ mapping $C^\theta_\ast[\mathscr S]$ to $C^{\theta'}_\ast[\mathscr S]$. Finally, set \eqq{
		\widecheck{C}_\ast[\mathscr S]&=C^o_\ast[\mathscr S]\oplus C^s_\ast[\mathscr S]\nonumber\\
		\check\partial=\begin{pmatrix}
		\partial^o_o&\partial^u_o\partial^s_u\\
		\partial^o_s&\partial^s_s+\partial^u_s\partial^s_u
		\end{pmatrix}&:\widecheck{C}_\ast[\mathscr S]\to \widecheck{C}_{\ast-1}[\mathscr S].\label{eqn:bdry-def}}
		
\begin{prop}
$\check\partial$ has degree $-1$, and $\check\partial^2=0$.
\end{prop}
\begin{proof}
To see that $\check\partial$ has degree $-1$, focus for simplicity on the matrix entry $\partial^o_o$ (the others are similar). $\partial^o_o$ sends chains $\sigma:\Delta\to B$ to chains $\ev_+:\sigma\times_B \breve M^+[\mathscr S](B,B')\to B'$, where $B,B'$ are interior critical points. Here the degree of $\sigma$ in $C^o_\ast[\mathscr S]$ is $\dim\Delta+\ind B$, whereas the degree of $\ev_+$ is \eq{
	\deg(\ev_+)&=\dim\left(\sigma\times_B \breve M^+[\mathscr S](B,B')\right)+\ind(B')\\
	&=\dim\Delta +\dim M^+[\mathscr S](B,B') - \dim B +\ind(B') \qquad\text{(by transversality)}\\
	&=\deg \sigma -1. \qquad\text{(by Proposition \ref{prop:delta-chain-with-dim})}}

The proof that $\check\partial^2=0$ is done in the gradient case in Proposition 4.2.4 in \cite{lin}. The proof carries over to the quasi-gradient setting, as we can substitute Proposition \ref{compactification} above to calculate the boundary of a moduli space.
\end{proof}
\begin{definition}
We define $\widecheck{H}_\ast[\mathscr S](X,v)$ to be the homology of the complex $(\widecheck{C}_\ast[\mathscr S],\check\partial)$ defined above.
\end{definition}

\subsection{Continuation maps}
\label{sec:continuation-maps}
Let $\{v_\lambda\}_{0\le \lambda\le 1}$ be a homotopy with $v_0, v_1$ Morse-Bott-Smale quasi-gradients. Suppose $N$ is a common isolating neighborhood for the flow of each $v_\lambda$, and set $\mathscr S_\lambda=\text{Inv}(N,v_\lambda)$. The purpose of this section is to prove the following proposition:
\begin{prop}
For generic homotopies $v_\lambda$, there is an isomorphism \eq{
	\widecheck{H}_\ast[\mathscr S_0](X,v_0)\simeq \widecheck{H}_\ast[\mathscr S_1](X,v_1).} 
\label{htopy-inv}
\end{prop}
\begin{cor}
The homology $\widecheck{H}_\ast[\mathscr S](X,v)$ is isomorphic to the reduced homology of the Conley index of $\mathscr S$ (with $\Z/2$ coefficients).
\end{cor}
\begin{proof}
Choose a small perturbation $v_1$ of $v$ which is Morse-Smale, and choose a generic homotopy $v_\lambda$ between $v$ and $v_1$. If $v_1$ is sufficiently close to $v$, then an isolating neighborhood $N$ for $\mathscr S$ remains an isolating neighborhood for each $v_\lambda$. According to the proposition, $\widecheck{H}_\ast[\mathscr S](X,v)\simeq \widecheck{H}_\ast[\mathscr S_1](X,v_1)$, where $\mathscr S_1=\text{Inv}(N,v_\lambda)$. Since $v_1$ is Morse-Smale, $\widecheck{H}_\ast[\mathscr S_1](X,v_1)$ is isomorphic to the reduced homology of the Conley index of $\mathscr S_1$ for the flow of $v_1$ [see \S 2.8 in \cite{lm}], but by homotopy invariance of the Conley index, this is isomorphic to the reduced homology of the Conley index of $\mathscr S_0$ for the flow of $v_0$.
\end{proof}

To prepare the proof of Proposition \ref{htopy-inv}, we first define a chain map \eq{
	\check F:\widecheck{C}_\ast[\mathscr S_0](X,v_0)\to \widecheck{C}_\ast[\mathscr S_1](X,v_1).} 

First, reparametrize $\lambda$ from $[0,1]$ to $[-\infty, \infty]$. Let \eqq{
	\rho:\R\to(-1,1) \label{increasing-diffeo}} be an increasing diffeomorphism, and define a vector field $\tilde v$ on $X\times \R$ by $\tilde v(x,\lambda)=(v_\lambda(x),\rho'(\lambda))$. Given stationary submanifolds $B_0$ and $B_1$ of $v_0$ and $v_1$, respectively, let $W^u(B_0), W^s(B_1)$ denote the set of flow lines of $\tilde v$ which limit to $B_0$ at $-\infty$ and $B_1$ at $+\infty$, respectively, and remain within $N$ for all time. For generic $\rho$, the intersections $W^u(B_0)\cap W^s(B_1)$ will be transverse for all pairs $(B_0,B_1)$, so we obtain smooth moduli spaces of unparametrized flow lines \eq{
	\breve M_\lambda(B_0,B_1)=\set{\gamma:\R\to X\times \R \mid \dot\gamma=\tilde v(\gamma), \lim_{t\to-\infty}\gamma(t)\in B_0,\;\lim_{t\to\infty}\gamma(t)\in B_1,\;\gamma(x,\lambda)\in N\text{ for all }x,\lambda}/\R}
which remain in $N$.

The following is analogous to Theorem 6.4 in \cite{bh}:
\begin{prop}
The moduli space $\breve M_\lambda(B_0, B_1)$ admits a compactification $\breve M_\lambda^+(B_0, B_1)$ consisting of broken trajectories from $B_0$ to $B_1$.
\label{second-compactification}
\end{prop}

As in Proposition 3.6.14 in \cite{lin}, we obtain
\begin{prop}
The compactified moduli space $\breve M_\lambda^+(B_0, B_1)$ defines an abstract $\delta$-chain of dimension $\ind(B_0)-\ind(B_1)$ (if $(B_0, B_1)$ is not boundary-obstructed) or $\ind(B_0)-\ind(B_1)+1$ (if $(B_0, B_1)$ is boundary-obstructed).
\label{interpolating-delta-chain}
\end{prop}

Elements of $\breve M_\lambda^+(B_0, B_1)$ can be written as broken trajectories $(\bm{\gamma_0}, \gamma, \bm{\gamma_1})$, where $\bm{\gamma_i}$ is a broken trajectory of the lifted vector field $\tilde v_i(x,\lambda)=(v_i(x),\rho'(\lambda))$ and $\gamma$ is a trajectory of $\tilde v$. The codimension-1 stratum of $\breve M_\lambda^+(B_0, B_1)$ admits the following description:

\begin{prop}
The codimension-1 stratum of $\breve M_\lambda^+(B_0, B_1)$ consists of broken trajectories $(\bm{\gamma_0}, \gamma, \bm{\gamma_1})$, where $\bm{\gamma_i}=(\gamma_i^1,\dots,\gamma_i^{n_i})$, for which there are exactly two trajectories in the set \eq{
    \set{\gamma_0^1,\dots, \gamma_0^{n_0}, \gamma, \gamma_1^1,\dots, \gamma_1^{n_1}}}
which are not boundary-obstructed.

\label{prop:codim-one-findim}
\end{prop}
\begin{proof}
This is analogous to Proposition 3.6.12 in \cite{lin}, or more specifically the special case of a cylindrical cobordism.
\end{proof}

Recall that $\widecheck{C}_\ast(X,v_0)$ is defined using $\delta$-chains transverse to the family \eq{
\mathcal F_{0}=\set{(\breve M^+(B,B'),\ev_-)\mid B,B'\in\mathfrak C[\mathscr S_0](v_0)}.} Enlarge $\mathcal F_0$ to $\hat{\mathcal F}$ by including the $\delta$-chains 
$(\breve M_\lambda^+(B_0, B_1),\ev_-)$ (with $B_0\in \mathfrak C[\mathscr S_0](v_0)$, $B_1\in \mathfrak C[\mathscr S_1](v_1)$).

The following lemma is a restatement of Lemma 4.2.12 in \cite{lin} in the current setting. 
\begin{lemma}
Write $\widecheck{C}^{\mathcal F}_\ast$ for the complex defined using the family $\mathcal F$. Then the inclusion $\widecheck{C}_\ast^{\hat{\mathcal F}}(X,v_0)\to \widecheck{C}_\ast^{\mathcal F_0}(X,v_0)$ is a quasi-isomorphism.
\label{enlarging-family}
\end{lemma}

For $\theta,\theta'\in \set{o,s,u}$ and a $\delta$-chain $\sigma$ in $B_0\in\mathfrak C^\theta[\mathscr S_0](v_0)$, we can define \eq{
	F^\theta_{\theta'}(\sigma)=\sum_{B_1\in\mathfrak C^{\theta'}[\mathscr S_1](v_1)} \sigma\times_{ B_0}\breve M_\lambda^+(B_0,B_1).}
We then set \eqq{
	\check F=\begin{pmatrix}
 F^o_o& F^u_o\partial^s_u+\partial^u_oF^s_u\\
 F^o_s& F^s_s+F^u_s\partial^s_u+\partial^u_sF^s_u \end{pmatrix}:\widecheck{C}_\ast[\mathscr S_0](X,v_0)\to \widecheck{C}_\ast[\mathscr S_1](X,v_1).\label{cont-map}}
 
 \begin{prop}
 The map $\check F$ defined above is a degree 0 chain map.
 \label{chain-map}
 \end{prop}
 \begin{proof}
 The proof is exactly analogous to Proposition 4.3.2 in \cite{lin}, and essentially follows from the description of the codimension-1 stratum of $M_\lambda^+(B_0,B_1)$ given in Proposition \ref{prop:codim-one-findim}.
 \end{proof}
 Note that the chain map $\check F$ has an implicit a priori dependence on the increasing function $\rho$ chosen in (\ref{increasing-diffeo}).
 
 To show that $\check F$ is an isomorphism in homology, it will be enough to show that
\begin{enumerate}[(1)]
\item if $v_0=v_1$ and $v_\lambda$ is constant in $\lambda$, then $\check F$ is the identity, and
\item if $v_2$ is a third Morse-Bott-Smale function and $\check F_{ij}:\widecheck{C}_\ast[\mathscr S_i](X,v_i)\to \widecheck{C}_\ast[\mathscr S_j](X,v_j)$, then $\check F_{12}\circ \check F_{01}$ is chain homotopic to $\check F_{02}$.
\end{enumerate}

\begin{prop}
If $v_0=v_1$ and we choose $v_\lambda$ to be constant in $\lambda$, then the chain map $\check F$ in (\ref{cont-map}) is the identity.
\label{const-id}
\end{prop}
\begin{proof}[Proof (cf. proof of Proposition 4.3.5 in \cite{lin})]
If $B,B'$ are two distinct critical submanifolds of $v_0$, then $\breve M_\lambda^+(B,B')$ (the collection of broken flow lines of $\tilde v$ from $B$ to $B'$) has a free $\R$-action given by translation in the $\lambda$-direction. Therefore, the evaluation maps $\ev_-:\breve M_\lambda^+(B,B')\to B'$ factor through the lower-dimensional space $\breve M_\lambda^+(B,B')/\R$. This implies that the $\delta$-chain $(\breve M_\lambda^+(B,B'),\ev_-)$ is negligible, in the sense of Definition \ref{crit-chains}, and hence that the moduli spaces $\breve M_\lambda^+(B,B')$ with $B,B'$ distinct do not contribute to $\check F$. On the remaining moduli spaces $\breve M_\lambda^+(B,B)$, the maps $\ev_{\pm}$ are homeomorphisms and so $\check F$ is the identity. 
\end{proof}

Now assume we are given a third Morse-Bott-Smale quasi-gradient $v_2$ as in item (2), and let $\rho_{ij}$ be the increasing function and $v_\lambda^{(ij)}$ the homotopy used in the construction of $\check F_{ij}$. Choose a smooth function $\rho_{012}:\R^2\to\R$ which is increasing in each coordinate and which satisfies \eq{
	\text{grad } \rho_{012}(\lambda,\mu)=\begin{cases}
	\rho_{01}'(\lambda)\partial_\lambda&\mu\ll 0\\
	\rho_{02}'(\mu)\partial_\mu&\lambda\ll 0\\
	\rho_{12}'(\mu)\partial_\mu&\lambda\gg 0\;. \end{cases}}
Choose a vector field $\hat v'$ on $X\times \R^2$ so that \eq{
	\hat v'(x,\lambda,\mu)&=\begin{cases}
	v^{(01)}_\lambda(x)&\mu\ll 0\\
	v^{(02)}_\mu(x)&\lambda\ll 0\\
	v^{(12)}_\mu(x)&\lambda \gg 0\\
	v_2(x)&\mu\gg 0\;.\end{cases}}
Now set $\hat v=\hat v'-\text{grad }\rho_{012}$.

Given stationary submanifolds $B_0$, $B_2$ of $v_0$ and $v_2$, let $W^u_{\lambda\mu}(B_0)$ be the set of flow lines $\gamma(a)$ of $\hat v$ limiting to $B_0\times\set{-\infty}\times\set{-\infty}$ as $a\to -\infty$, and let $W^s_{\lambda\mu}(B_2)$ be the set of flow lines $\gamma(a)$ of $\hat v$ limiting to $B_2\times\set{\infty}\times\set{\infty}$ as $a\to \infty$. For a generic choice of $\rho_{012}$, the intersection $W^u_{\lambda\mu}(B_0)\cap W^s_{\lambda\mu}(B_2)$ is transverse and we obtain a smooth moduli space $\breve M_{\lambda\mu}(B_0, B_2)=(W^u_{\lambda\mu}(B_0)\cap W^s_{\lambda\mu}(B_2))/\R$.

As with Propositions \ref{second-compactification} and \ref{interpolating-delta-chain} above, the following two propositions follow as in Theorem 6.4 in \cite{bh} and Proposition 3.6.14 in \cite{lin}.
\begin{prop}
The moduli space $\breve M_{\lambda\mu}(B_0, B_2)$ admits a compactification $\breve M_{\lambda\mu}^+(B_0, B_2)$ consisting of 
\begin{enumerate}[(1)]
\item broken trajectories of $\hat v$ from $B_0$ to $B_2$;
\item broken trajectories of $\tilde v^{(02)}$ from $B_0$ to $B_2$;
\item concatenations of broken trajectories of $\tilde v^{(01)}$ from $B_0$ to an intermediate stationary submanifold $B_1$ of $v_1$ with broken trajectories of $\tilde v^{(12)}$ from $B_1$ to $B_2$.
\end{enumerate}
\end{prop}

\begin{prop}
The compactified moduli space $\breve M_{\lambda\mu}^+(B_0, B_2)$ defines an abstract $\delta$-chain.
\label{chain-homotopy-delta}
\end{prop}

As before, for $\theta,\theta'\in \set{o,s,u}$ and a $\delta$-chain $\sigma$ in $B_0\in\mathfrak C^\theta[\mathscr S_0](v_0)$, we can now define \eq{
	\Psi^\theta_{\theta'}(\sigma)&=\sum_{B_2\in\mathfrak C^{\theta'}[\mathscr S_2](v_2)} \sigma\times_{ B_0}\breve M_{\lambda\mu}^+(B_0,B_2)\\
	\check \Psi&=\begin{pmatrix}
 \Psi^o_o& \Psi^u_o\partial^s_u+\partial^u_o\Psi^s_u\\
 \Psi^o_s& \Psi^s_s+\Psi^u_s\partial^s_u+\partial^u_s\Psi^s_u \end{pmatrix}:\widecheck{C}_\ast[\mathscr S_0](X,v_0)\to\widecheck{C}_\ast[\mathscr S_2](X,v_2).}

\begin{prop}
The map $\Psi$ satisfies $\partial\Psi+\Psi\partial=\check F_{fh}+\check F_{gh}\circ\check F_{fg}$.
\label{cont-htopy-map}
\end{prop}
\begin{proof}
The proof is entirely analogous to Proposition 3.6 in \cite{lin}.
\end{proof}
\begin{proof}[Proof of Proposition \ref{htopy-inv}]
Let $v_\lambda$ be a homotopy between $v_0$ and $v_1$. Let $v'_\lambda$ be the reversed homotopy from $v_1$ to $v_0$, and let $v_\lambda^0$ be the constant homotopy at $v_0$. By Proposition \ref{chain-map} these give rise to chain maps $\check F_{01}$, $\check F_{10}$, and $\check F_{00}$, respectively, of which the last is the identity map by Proposition \ref{const-id}. By Proposition \ref{cont-htopy-map}, $\check F_{10}\circ \check F_{01}$ is chain homotopic to $\check F_{00}$, which means $\check F_{10}$, $\check F_{01}$ are inverses in homology.
\end{proof}

\subsection{Morse-Bott homology for manifolds with $\pin(2)$-actions}\label{subsection:mb-homology-mfld-pin2}
\subsubsection{Real-oriented blowups}
Let $X$ now be a closed Riemannian manifold with an action of the group \eq{
	G=\pin(2)=S^1\cup j\cdot S^1\subset\mathbb H.} by isometries. We assume this action is \textit{semifree}, i.e. the stabilizer of a point in $X$ is $\pin(2)$ or $\{1\}$. Let $\Fps$ be the fixed set of the action. As in \S 2.6 of \cite{lm}, we can construct the real-oriented blowup \eq{
	X^\sigma=(X-\Fps)\cup \left(N^1(\Fps)\times[0,\epsilon)\right),} where $N^1(\Fps)$ is the unit normal bundle of $\Fps$ in $X$. The action of $\pin(2)$ on $X^\sigma$ is free by construction. 
	
	For $\fps\in \Fps$, write $L_\fps=(\nabla\tilde v)|_{N_\fps\Fps}:N_\fps\Fps\to N_\fps\Fps$. Note that $N_\fps\Fps$ naturally has the structure of a vector space over $\mathbb H$, and that $L_\fps$ is $\mathbb H$-linear. We define a Morse-Bott equivariant quasi-gradient in our situation as in Definition 2.18 in \cite{lm}:
\begin{definition}
A smooth, $\pin(2)$-equivariant vector field $\tilde v$ on $X$ is called a \textit{Morse-Bott equivariant quasi-gradient} if the following conditions are satisfied:
\begin{enumerate}[(1)]
\item All stationary points of the induced vector field $v$ on $(X-\Fps)/S^1$ are hyperbolic.
\item All stationary points of $\tilde v|_\Fps$ are hyperbolic.
\item At each stationary point $\fps$ of $\tilde v|_\Fps$, $L_\fps$ is self-adjoint and admits a basis over $\mathbb H$ of eigenvectors $\phi_1(\fps),\dots,\phi_n(\fps)$ with corresponding eigenvalues $\lambda_1(\fps),\dots, \lambda_n(\fps)$ such that \eq{
	\lambda_1(\fps)<\dots<\lambda_n(\fps)}
and $\lambda_i(\fps)\ne 0$ for any $i$.
\item There exists a smooth, $\pin(2)$-equivariant function $\tilde f:X\to \R$ such that $d\tilde f(\tilde v)\ge 0$ at all $x\in X$, with equality holding if and only if $\tilde v(x)=0$.
\end{enumerate}
Such a vector field is called \textit{Morse-Bott-Smale} if the induced vector field $v^\sigma$ on $X^\sigma/S^1$ satisfies the Morse-Bott-Smale condition for boundary-unobstructed trajectories and the Morse-Bott-Smale condition in $\partial(X^\sigma/S^1)$ for boundary-obstructed trajectories.
\label{morse-eq-qg}
\end{definition}

\begin{remark}
The stationary points of the induced vector field $ v^\sigma$ on $\partial(X^\sigma)$ correspond to pairs $(\fps,\phi)$, where $\fps\in \Fps$ and $\phi$ is a unit-length eigenvector of $L_\fps$. Since $L_\fps$ is $\mathbb H$-linear, there is an $S^3$'s worth of unit eigenvectors for each eigenvalue, which means the critical loci of $v^\sigma$ on $X^\sigma/S^1$ are isolated points (in the interior) and copies of $S^3/S^1=S^2$ (on the boundary).
\end{remark}
\begin{lemma}
Parts (1), (2), and (3) of Definition \ref{morse-eq-qg} are equivalent to asking that the induced vector field $v^\sigma$ on $X^\sigma/S^1$ is hyperbolic in normal directions, i.e. for every stationary submanifold $B\subset X^\sigma/S^1$ of $v^\sigma$ and for every $x\in B$, the operator \eq{
    L_x^\sigma=(\nabla v^\sigma)|_{N_x(B)}:N_x(B)\to N_x(B)} is hyperbolic.
\label{mb-eq-qg-lemma}
\end{lemma}
\begin{proof}
Clearly condition (1) is equivalent to interior stationary points of $v^\sigma$ being hyperbolic. Therefore we restrict our attention to the boundary and fix a boundary stationary point $x=(r,[\phi])$ of $v^\sigma$, where $L_r\phi=\lambda \phi$. 
Following Lemma 2.5.5 of \cite{m3m}, we may decompose the tangent space to $X^\sigma/S^1$ at $x$ as \eq{
    T_x(X^\sigma/S^1) &=T_\fps\Fps\oplus (\C \phi)^\perp \oplus \R}
and write the linearization of $v^\sigma$ as \eq{
    \nabla v^\sigma(x)&=\begin{bmatrix}
    (\nabla \tilde v)|_R & 0 & 0\\
    \ast & L_r-\lambda & 0\\
    0&0& \lambda\end{bmatrix}.}
The stationary submanifold $B$ containing $x$ has $T_xB=\set{0}\oplus \C (j\phi)\oplus \set{0}$, so $N_xB=T_rR\oplus(\C\phi)^\perp/\C(j\phi)\oplus\R$ and \eq{
    L_x^\sigma= \begin{bmatrix}
    (\nabla \tilde v)|_R & 0 & 0\\
    \ast & L_r-\lambda & 0\\
    0&0& \lambda\end{bmatrix}:N_xB\to N_xB.}
This is hyperbolic if and only if $\tilde v|_R$ is hyperbolic, $\lambda\ne 0$, and $\lambda$ has multiplicity exactly 2, so that $L_r-\lambda$ is invertible on the complement of $\C\phi\oplus\C(j\phi)$.

\end{proof}

\begin{lemma}
Let $\tilde v$ be a Morse-Bott equivariant quasi-gradient vector field on $X$, as in Definition \ref{morse-eq-qg}. Let $\gamma:\R\to X$ be a flow line of $\tilde v$. Then $\lim_{t\to\infty}[\gamma(t)]$ and $\lim_{t\to -\infty}[\gamma(t)]$ exist in $X/S^1$, and they are both projections of stationary points of $\tilde v$.
\label{limits-exist-in-quot}
\end{lemma}
\begin{proof}
This is identical to Lemma 2.22 in \cite{lm}, which is stated for Morse equivariant quasi-gradients. In the current setting, this lemma is proved analogously to Lemma \ref{limits-exist-noneq}.
\end{proof}
\begin{lemma}
Let $\tilde v$ be a Morse-Bott equivariant quasi-gradient vector field on $X$, as in Definition \ref{morse-eq-qg}. Let $\gamma:\R\to X^\sigma/S^1$ be a flow line of $v^\sigma$. Then $\lim_{t\to\infty}\gamma(t)$ and $\lim_{t\to -\infty}\gamma(t)$ exist in $X^\sigma/S^1$, and they are both stationary points of $v^\sigma$.
\label{limits-exist-eq}
\end{lemma}

\begin{proof}
This is similar to the proof of Lemma 2.23 in \cite{lm}, which is reproduced below with the necessary changes.

It suffices to consider the half-trajectory $\gamma_+=\gamma|_{\R_+}$. By Lemma \ref{limits-exist-in-quot}, the projection of $\gamma_+$ to the blow-down $X/S^1$ limits to a stationary point $r_0$ at $\infty$. Since the projection is one-to-one away from the boundary, we can assume $r_0$ is in the fixed-point set $R$ of the $\pin(2)$-action. 

A small neighborhood $V\subset X^\sigma/S^1$ of the boundary $\partial(X^\sigma/S^1)$ can be identified with the normal bundle to the boundary. Hence a point $v\in V$ can be written as a triple $(r,s,[\phi])$ with $r\in R$, $s\ge 0$, and $\phi\in N_r^1(R)$, normalized so that $|\phi|=1$. Define \eq{
    \Lambda:V\to \R,\quad \Lambda(r,s,[\phi])=\inner{\phi,L_r\phi}.}
The restriction of $v^\sigma$ to $F=N_{r_0}^1(R)/S^1$ is the gradient of $\frac12\Lambda|_{F}$ (see e.g. the discussion above Definition 2.18 in \cite{lm}). Therefore, $v^\sigma|_F$ has $n$ stationary submanifolds, namely 2-spheres labeled by the eigenvalues $\lambda_i$ of $L_r$.

Since the blow-down of $\gamma_+$ converges to $r_0$, given any neighborhood $U$ of $F$, there is $t_0$ such that $\gamma_+(t)\in U$ for all $t>t_0$. In particular, all the accumulation points of $\gamma_+$ as $t\to\infty$ must be contained in $F$. It remains to show that there is a unique accumulation point, and that that point is a stationary point of $v^\sigma$.

If $z\in F$ is an accumulation point of $\gamma_+$, then so are all the points on the flow trajectory $\zeta$ through $z$. Now $\zeta$ is contained in $F$, where the flow is a gradient flow, and hence $\zeta$ limits to stationary points $x,y$ at $-\infty$ and $+\infty$, respectively, corresponding to eigenvectors $[\phi]$ and $[\psi]$ of $L_{r_0}$. If $[\phi]\ne [\psi]$, then the eigenvalues $\lambda_i, \lambda_j$ corresponding to $[\phi],[\psi]$ must also be distinct, because trajectories of $\Lambda$ do not connect distinct eigenvectors with the same eigenvalue. Note that $x,y$ are also accumulation points for $\gamma_+$, and $x=y$ if and only if $\zeta$ is a stationary trajectory.

If $x\ne y$, and hence $\lambda_i\ne \lambda_j$, pick $\lambda$ with $\lambda_i > \lambda >\lambda_j$ such that $c$ is not an eigenvalue of $L_r$. Since $v^\sigma|_F$ is the gradient of $\frac12\Lambda$, we have \eq{
    d\Lambda(v^\sigma)>0\;\text{ on } \Lambda^{-1}(\lambda)\cap F.}
By compactness of $\Lambda^{-1}(\lambda)\cap F$, we can find a neighborhood $U\subset V$ of $F$ such that \eqq{
    d\Lambda(v^\sigma)>0\;\text{ on } \Lambda^{-1}(\lambda)\cap U.\label{Lambda-decr}}
Since $x,y$ are accumulation points of $\gamma_+$, the half-trajectory $\gamma_+$ must intersect $\Lambda^{-1}(\lambda)$ infinitely many times. Furthermore, for $t\ge t_0$ the trajectory $\gamma_+$ is contained in $U$. This contradicts equation (\ref{Lambda-decr}), which says that $\Lambda$ decreases every time the trajectory goes through $\Lambda^{-1}(\lambda)\cap U$.
\end{proof}
As in \S 2.6 of \cite{lm}, the assumption that $v$ is a Morse-Bott equivariant quasi-gradient means that even though the vector field $v^\sigma$ on $X^\sigma/S^1$ is not a quasi-gradient vector field in any natural way, we can still use it to compute the Morse-Bott complex $\widecheck C_\ast(X^\sigma/S^1, v^\sigma)$: Lemma \ref{mb-eq-qg-lemma} shows that the stationary manifolds of $v^\sigma$ are normally hyperbolic, which is enough to establish Proposition \ref{compactification}.

The complex $\widecheck C_\ast(X^\sigma/S^1,v^\sigma)$ computes the singular homology of $X^\sigma/S^1$. However, it still has a residual $\mathbb Z/2$-action by $j$. 
\begin{definition}
We define \eq{
    \widecheck C_\ast^{\text{inv}}(X^\sigma/S^1,v^\sigma)\subset C_\ast(X^\sigma/S^1,v^\sigma)} to be the subcomplex of $j$-invariant chains, and denote \eq{
    \widecheck H_\ast(X,\tilde v)=H_\ast(C_\ast^{\text{inv}}(X^\sigma/S^1,v^\sigma)).}
If $\scr S\subset X$ is a $\pin(2)$-invariant isolated invariant set for the flow of $\tilde v$, then $\scr S^\sigma/S^1\subset X^\sigma/S^1$ is an isolated invariant set for the flow of $v^\sigma$, and we define \eq{
    \widecheck H_\ast[\scr S](X,\tilde v)=H_\ast(C_\ast^{\text{inv}}[\scr S^\sigma/S^1](X^\sigma/S^1,v^\sigma)).}
\end{definition}
\begin{remark}
There is a clear quasi-isomorphism between the complex $\widecheck C_\ast^{\text{inv}}(X^\sigma/S^1,v^\sigma)$ and the complex $\widecheck C_\ast(X^\sigma/\pin(2),v^\sigma)$ obtained by considering stationary points and trajectories in the full quotient $X^\sigma/\pin(2)$. (We generally prefer to work in the intermediate quotient $X^\sigma/S^1$, however.) Hence, \eq{
    \widecheck H_\ast(X,\tilde v)\simeq\tilde H_\ast(X^\sigma/\pin(2))\simeq\tilde H_\ast^{\pin(2)}(X^\sigma),} and thus $\widecheck H_\ast(X,\tilde v)$ approximates the $\pin(2)$-equivariant homology of $X$. More precisely, since $X^\sigma$ is $\pin(2)$-equivariantly homotopy equivalent to $X-X^{S^1}$, we see that \eq{
    \tilde H_\ast^{\pin(2)}(X^\sigma)\simeq\tilde H_\ast^{\pin(2)}(X-X^{S^1}),} and thus \eq{
    \widecheck H_j(X,\tilde v)\simeq \tilde H_j^{\pin(2)}(X)\quad (j\le n-1),} where $n$ is the connectivity of the pair $(X,X-X^{S^1})$.
\label{rmk:no-inv-set}
\end{remark}

\begin{remark}
To develop Remark \ref{rmk:no-inv-set} further, suppose $\mathscr S\subset X$ is an isolated invariant set of the flow of $\tilde v$ which is $\pin(2)$-invariant. According to \cite{pruszkoConleyIndexFlows1999} and \cite{floerRefinementConleyIndex1987}, we may choose a $\pin(2)$-invariant index pair $(N,L)$ for $\mathscr S$. Then $\mathscr S^\sigma/\pin(2)$ is an isolated invariant set for the vector field on $X^\sigma/\pin(2)$ induced by $v^\sigma$, and hence has a Conley index $I(\mathscr S^\sigma/\pin(2))$. We then have \eq{
    \widecheck H_\ast[\mathscr S](X,\tilde v)\simeq \tilde H_\ast(I(\mathscr S^\sigma/\pin(2))).}
But $(N^\sigma/\pin(2),L^\sigma/\pin(2))$ is an index pair for $\scr S^\sigma/\pin(2)$, and so \eq{
    \tilde H_\ast(I(\mathscr S^\sigma/\pin(2)))\simeq\tilde H_\ast^{\pin(2)}(N^\sigma/L^\sigma).} 
    
Since the blowup $N^\sigma$ is $\pin(2)$-equivariantly homotopy equivalent to $N-N^{S^1}$, and similarly $L^\sigma \simeq L-L^{S^1}$, we see that \eq{
    \widecheck H_\ast[\mathscr S](X,\tilde v)\simeq\tilde H_\ast^{\pin(2)}(I(\scr S)-I(\scr S)^{S^1}),} and so \eq{
    \widecheck H_j[\mathscr S](X,\tilde v)\simeq\tilde H_j^{\pin(2)}(I(\scr S)) \quad (j\le n-1),} where $n$ is the connectivity of the pair $(I(\scr S), I(\scr S)-I(\scr S)^{S^1})$.
\label{rmk:computes-eqhom-spec}
\end{remark}

\subsubsection{The $Q$- and $V$-actions in Morse-Bott homology}\label{subsubsection:QV} Recall that \eq{
	H^\ast(B\pin(2);\mathbb F)=\mathbb F[Q,V]/(Q^3),} with $\deg Q=1$, $\deg V=4$. With $(X,\tilde v)$ as in the previous subsection, we now define chain maps of degree $-1$ and $-4$ on $\widecheck C_\ast^{\text{inv}}(X^\sigma/S^1,\tilde v)$ corresponding to the $H^\ast(B\pin(2);\mathbb F)$-module structure on $H_\ast^{\pin(2)}(X^\sigma)$. The $\pin(2)$-action on $X^\sigma$ produces a natural quaternionic line bundle $E^\sigma$ on $X^\sigma/\pin(2)$. The $\Z/2$-bundle $\pi:X^\sigma/S^1\to X^\sigma/\pin(2)$ also gives rise to a natural real line bundle $F^\sigma$ over $X^\sigma/\pin(2)$. Let $\zeta,\eta$ be smooth sections of $E^\sigma$ and $F^\sigma$ (respectively) which are transverse to the zero sections of $E^\sigma$ and $F^\sigma$ both in $X^\sigma/\pin(2)$ and within $\partial(X^\sigma/\pin(2))$. Let $Z_\zeta$, $Z_\eta$ denote the zero sets of $\zeta,\eta$. Then for stationary submanifolds $B$, $B'$ of $v^\sigma$, let $M(B,B')$ denote the usual moduli space of flow lines $\gamma:\R\to X^\sigma/S^1$ from $B$ to $B'$, and define cut-down moduli spaces \eq{
	M^{\cap\zeta}(B,B')&=\{\gamma\in M(B,B')\mid \pi(\gamma(0))\in Z_\zeta\}\\
	M^{\cap\eta}(B,B')&=\{\gamma\in M(B,B')\mid  \pi(\gamma(0))\in Z_\eta\}.} The cut-down moduli spaces define abstract $\delta$-chains. By Lemma \ref{enlarging-family}, we may enlarge the family $\mathscr F$ to include the $\delta$-chains $(M^{\cap\zeta}(B,B'),\ev_-),(M^{\cap\eta}(B,B'),\ev_-)$ for all $B,B'$. Denote by $\mathfrak C^o(X^\sigma/S^1$, $\mathfrak C^s(X^\sigma/S^1)$, and $\mathfrak C^u(X^\sigma/S^1)$ the collection of interior, boundary-stable, and boundary-unstable stationary submanifolds for $v^\sigma$, respectively. With notation as in equation (\ref{breakup}), and for $\theta,\theta'\in\set{o,s,u}$ and a $\delta$-chain $\sigma$ in $B\in\mathfrak C^\theta(X^\sigma/S^1)$, define \eq{
	m(\zeta)^\theta_{\theta'}: C^\theta_\ast(X^\sigma/S^1)\to C^{\theta'}_{\ast-4}(X^\sigma/S^1)\\
	 m(\eta)^\theta_{\theta'}: C^\theta_\ast(X^\sigma/S^1)\to C^{\theta'}_{\ast-1}(X^\sigma/S^1)} by \eq{
	m(\zeta)^\theta_{\theta'}(\sigma)=\sum_{B'\in\mathfrak C^\theta(X^\sigma/S^1)} \sigma\times_{\ev_-} M^{\cap\zeta}(B,B')\\
	m(\eta)^\theta_{\theta'}(\sigma)=\sum_{B'\in\mathfrak C^\theta(X^\sigma/S^1)} \sigma\times_{\ev_-} M^{\cap\eta}(B,B'),} where each summand is viewed as a $\delta$-chain in $B'$ via $\ev_+$. Finally, restrict to invariant chains and define \eq{
	\check m(\zeta)&=\begin{pmatrix} m(\zeta)^0_0 & -m(\zeta)^u_0\partial_u^s-\partial_o^u m(\zeta)^s_u\\
	m(\zeta)^0_s & m(\zeta)^s_s -m(\zeta)^u_s\partial^s_u-\partial^u_s m(\zeta)^s_u \end{pmatrix}:\widecheck{C}^{\text{inv}}_\ast(X^\sigma/S^1)\to \widecheck{C}^{\text{inv}}_{\ast-4}(X^\sigma/S^1)\\
	\check m(\eta)&=\begin{pmatrix} m(\eta)^0_0 & -m(\eta)^u_0\partial_u^s-\partial_o^u m(\eta)^s_u\\
	m(\eta)^0_s & m(\eta)^s_s -m(\eta)^u_s\partial^s_u-\partial^u_s m(\eta)^s_u \end{pmatrix}:\widecheck{C}^{\text{inv}}_\ast(X^\sigma/S^1)\to\widecheck{C}^{\text{inv}}_{\ast-1}(X^\sigma/S^1).}
The maps induced by $\check m(\zeta)$, $\check m(\eta)$ on the homology $\widecheck H_\ast(X,\tilde v)$ are exactly the cap products with $p_1(E^\sigma)$ and $w_1(F^\sigma)$, because $Z_\zeta$ (respectively $Z_\eta$) is Poincar\'e dual to $p_1(E^\sigma)$ (resp. $w_1(F^\sigma)$), which is the pullback of $V$ (resp. $Q$) from $H^\ast (B\pin(2);\bb F)$ to $H^\ast(X^\sigma/\pin(2);\bb F)$. This recovers the $\mathbb F[Q,V]/(Q^3)$-module structure on $H_\ast^{\pin(2)}(X^\sigma)$.

The isomorphism in Remark \ref{rmk:computes-eqhom-spec} respects the actions of $Q$ and $V$ so defined, leaving us with the following proposition:
\begin{prop}
Let $X$ be a finite-dimensional Riemannian manifold without boundary, equipped with a semifree action of the group $\pin(2)$. Suppose $\tilde v$ is a Morse-Bott equivariant quasi-gradient on $X$. Then for any $\pin(2)$-invariant isolated invariant set $\mathscr S$ for the flow of $\tilde v$, the $\pin(2)$-equivariant Conley index $I(\mathscr S)$ satisfies \eq{
    \widecheck H_j[\mathscr S](X,\tilde v)\simeq\tilde H_j^{\pin(2)}(I(\scr S)) \quad (j\le n-1),} where $n$ is the connectivity of the pair $(I(\scr S), I(\scr S)-I(\scr S)^{S^1})$. The isomorphism respects the actions of the ring $\mathcal R=\mathbb F[Q,V]/(Q^3)$ on either side.
    
\label{prop:fin-dim-summary}
\end{prop} 

\section{Seiberg-Witten preliminaries}\label{sec:SW-prelim}
This section briefly sets up the configuration spaces and operators relevant to the proof of Theorem \ref{thm:main-thm}. Throughout the section, we fix a rational homology 3-sphere $Y$ with a metric $g$ and a spin structure $\mathfrak s$ with associated spinor bundle $\mathbb S$.
\subsection{The blown-up configuration space}
    The affine space of spin$^c$ connections on $Y$ can be identified with $\Omega^1(Y; i\R)$, where the spin connection acts as the zero. The \textit{configuration space} $\mathcal C(Y)$ is then defined to be \eq{
    \mathcal C(Y)=\Omega^1(Y;i\R)\oplus\Gamma(\mathbb S),} and has an action of the gauge group \eq{
    \mathcal G(Y)=C^\infty(Y,S^1)} by $u\cdot(a,\phi)=(a-u^{-1}du,u\phi)$. 
    The blown-up configuration space $\mathcal C^\sigma(Y)$ is defined to be \eq{
    \mathcal C^\sigma(Y)=\set{(a,s,\phi)\mid s\ge 0,\|\phi\|_{L^2(\mathbb S)}=1}\subset \Omega^1(Y; i\R)\times \R_{\ge 0}\times \Gamma(\mathbb S),} and admits a blow-down map $\mathcal C^\sigma(Y)\to \mathcal C(Y)$, $(a,s,\phi)\mapsto (a,s\phi)$, which is a diffeomorphism off the \textit{reducible locus} $s=0$. The action of $\mathcal G(Y)$ on $\mathcal C(Y)$ induces a free action of $\mathcal G(Y)$ on $\mathcal C^\sigma(Y)$, and the quotient is \eq{
    \mathcal B^\sigma(Y)=\mathcal C^\sigma(Y)/\mathcal G(Y).}
    
    The spinor bundle $\mathbb S$ admits a right action by the quaternions $\bb H$, and in particular by $j\in\bb H$. Thus, we get a map \eq{
        \jmath: \mathcal C(Y)&\to\mathcal C(Y)\\
        (a,\phi)&\mapsto (-a,\phi\cdot j).} The map $\jmath$ extends naturally to $\mathcal C^\sigma(Y)$ and descends to an involution on $\mathcal B^\sigma(Y)$. 

\subsection{The perturbed Seiberg-Witten flow}
    The \textit{Seiberg-Witten flow} on $\mathcal C(Y)$ is given by the equations \eq{
    \frac{da}{dt}&=-\ast da-\tau(\phi,\phi)\\
    \frac{d\phi}{dt}&=-D_a\phi,}
    where $\tau$ is a quadratic term and $D_a$ is the Dirac operator for the connection corresponding to $a$. The \textit{Seiberg-Witten vector field} $\mathcal X$ is the vector field generated by this flow.

    The Seiberg-Witten flow on $\mathcal C(Y)$ induces a flow on the blow-up $\mathcal C^\sigma(Y)$, given by \eq{
        \frac{da}{dt}&=-\ast da-s^2\tau(\phi,\phi)\\
        \frac{ds}{dt}&=-\Lambda(a,s,\phi)s\\
        \frac{d\phi}{dt}&=-D_a\phi+\Lambda(a,s,\phi)\phi,}
    where \eq{
        \Lambda(a,s,\phi)=\text{Re}\inner{\phi,D_a\phi}_{L^2}.}
    This flow generates a vector field $\mathcal X^\sigma$. The Seiberg-Witten flow is invariant under both the action of $\mathcal G(Y)$ and the action of $\jmath$.

    Fix a regular $\jmath$-invariant perturbation \eq{
    \frak q:\mathcal C(Y)\to \mathcal C(Y),} obtained as the gradient of a $\pin(2)$-invariant cylinder function (see \cite{lin}, Definition 5.2.1 and Theorem 5.2.6). Write $\frak q^0$ and $\frak q^1$ for the 1-form and spinor components of $\frak q$, respectively. Then there is a perturbed Seiberg-Witten flow on $\mathcal C(Y)$ given by \eq{
    \frac{da}{dt}&=-\ast da-\tau(\phi,\phi)-\frak q^0(a,\phi)\\
    \frac{d\phi}{dt}&=-D_a\phi-\frak q^1(a,\phi),} and an induced flow on $\mathcal C^\sigma(Y)$ given by \eq{
        \frac{da}{dt}&=-\ast da-s^2\tau(\phi,\phi)-\frak q^0(a,s\phi)\\
        \frac{ds}{dt}&=-\Lambda_{\frak q}(a,s,\phi)s\\
        \frac{d\phi}{dt}&=-D_a\phi-\tilde{\frak q}^1(a,s\phi)+\Lambda_{\frak q}(a,s,\phi)\phi,}
    where \eq{
        \tilde{\frak q}^1(a,s,\phi)=\int_0^1\mathcal D_{(a,st\phi)}\frak q^1(0,\phi)\,dt} and \eq{
        \Lambda_{\frak q}(a,s,\phi)=\text{Re}\inner{\phi,D_a\phi+\tilde{\frak q}^1(a,s,\phi)}_{L^2}.} These flows give rise to vector fields $\mathcal X_{\frak q}$ and $\mathcal X^\sigma_{\frak q}$, both of which are still $\mathcal G(Y)$- and $\jmath$-invariant.

\subsection{Coulomb slices}
The global Coulomb slice in $\mathcal C(Y)$ is \eq{
    W=\ker d^\ast \oplus \Gamma(\bb S).} Any $(a,\phi)\in \mathcal C(Y)$ can be moved into the global Coulomb slice by a unique element of the \textit{normalized gauge group} \eq{
    \mathcal G^\circ(Y)=\set{e^f\Bigm\vert f:Y\to i\R,\int_Y f=0}\subset\mathcal G(Y).} The global Coulomb projection is given by \eq{
    \Pi^{\text{gC}}(a,\phi)=(a-df,e^f\phi)} where $f=Gd^\ast a$, $G$ being the Green's operator of $\Delta=d^\ast d$. There is still a residual gauge action on $W$ by the constant gauge transformations, and $\jmath$ preserves $W$.

    The \textit{enlarged local Coulomb slice} $\mathcal K^{\text{e}}_{(a,\phi)}$ at $(a,\phi)\in\mathcal C(Y)$ is the $L^2$-orthogonal complement to the orbits of the normalized gauge group $\mathcal G^\circ(Y)$. There is a well-defined \textit{enlarged local Coulomb projection} $\Pi^{\elc}_{(a,\phi)}:T_{(a,\phi)}\mathcal C(Y)\to\mathcal K^{\text{e}}_{(a,\phi)}$ (see \cite{lm}, Lemma 3.2), which is an isomorphism when restricted to $T_{(a,\phi)} W$. By measuring the $L^2$-inner product of the enlarged local Coulomb projections of tangent vectors to $W$, we obtain a metric $\tilde g$ on $W$. Finally, the \textit{anticircular global Coulomb slice} $\mathcal K^{\agc}_{(a,\phi)}$ at $(a,\phi)\in W$ is defined to be the $\tilde g$-orthogonal complement of the tangent to the $S^1$-orbit at $(a,\phi)$; explicitly, \eq{
    \mathcal K^{\agc}_{(a,\phi)}=\set{(b,\psi)\in T_{(a,\phi)}\mathcal C(Y)\Bigm\vert d^\ast b=0,\,\inner{(0,i\phi),(b,\psi)}_{\tilde g}=0}.}

    The global Coulomb slice $W$ has a blow-up $W^\sigma\subset \mathcal C^\sigma(Y)$, and $W^\sigma/S^1$ is a model for $\mathcal B^\sigma(Y)$. There is a vector field $\mathcal X_{\frak q}^{\gc}$ on $W$ obtained from $\mathcal X_{\frak q}$ by global Coulomb projection, and $\mathcal X_{\frak q}^{\gc}$ induces a vector field $\mathcal X_{\frak q}^{\gc,\sigma}$ on $W^\sigma$.

    The anticircular global Coulomb slices also extend naturally to slices $\mathcal K^{\agc,\sigma}$ in the blow-up forming an $S^1$-invariant bundle over $W^\sigma$. The slice $\mathcal K^{\agc,\sigma}_{(a,s,\phi)}$ is a model for the tangent space to $W^\sigma/S^1$ at the orbit of $(a,s,\phi)$. Since $\mathcal X_{\frak q}^{\gc,\sigma}$ is a section of $\mathcal K^{\agc,\sigma}$, it descends to a vector field $\mathcal  X_{\frak q}^{\agc,\sigma}$ on $W^\sigma/S^1$.

\subsection{Eigenvalue cutoffs}
The vector field $\mathcal X_{\frak q}^{\gc}$ has the form \eq{
    \mathcal X_{\frak q}^{\gc}=l+c_{\frak q},} where $l$ is a linear Fredholm operator, self-adjoint with respect to the $L^2$-metric, and $c_{\frak q}$ is a quadratic perturbation. For each $\lambda>1$, let $W^\lambda$ be the finite-dimensional subspace of $W$ spanned by eigenvectors of $l$ with eigenvalues in $(-\lambda,\lambda)$. Let $\tilde p^\lambda:W\to W^\lambda$ denote the $L^2$-orthogonal projection onto $W^\lambda$. As in \S3.4 of \cite{lm}, it is necessary to modify the $\tilde p^\lambda$ to produce maps $p^\lambda:W\to W^\lambda$ with the following properties:
    \begin{itemize}
    \item the $p^\lambda$ vary smoothly in $\lambda$;
    \item there is a sequence $\lambda_1^\bullet<\lambda_2^\bullet<\dots$ converging to infinity such that $p^{\lambda_i^\bullet}=\tilde p^{\lambda_i^\bullet}$.
    \end{itemize}
    This is accomplished as follows. First, modify the $\tilde p^\lambda$ to be smooth in $\lambda$ by defining \eq{
        p^\lambda_{\text{prel}}=\int_0^1\beta(\theta)\tilde p^{\lambda-\theta}\,d\theta,}
    where $\beta$ is a smooth, non-negative function which is nonzero exactly on $(0,1)$ and satisfies $\int_{\R}\beta(\theta)\,d\theta=1$. Next, fix a sequence $\lambda_1^\bullet<\lambda_2^\bullet<\dots$ converging to infinity
    such that no $\lambda_i^\bullet$ is an eigenvalue of $l$. Fix also disjoint intervals $[\lambda_i^\bullet-\epsilon_i,\lambda_i^\bullet+\epsilon_i]$ that do not contain eigenvalues of l, and let $\beta_i:(0,\infty)\to[0,1]$ be a smooth bump function supported in $[\lambda_i^\bullet-\epsilon_i,\lambda_i^\bullet+\epsilon_i]$ and satisfying $\beta_i(\lambda_i^\bullet)=1$. Finally, set \eq{
        p^\lambda&=\sum_i \beta_i(\lambda)\tilde p^\lambda+\left(1-\sum_i\beta_i(\lambda)\right)p^\lambda_{\text{prel}}.} It is clear that these $p^\lambda$ satisfy the two desired properties.
        
    We define an approximate Seiberg-Witten vector field on $W$ by \eq{
    \mathcal X_{\frak q^\lambda}^{\gc}=l+p^\lambda c_{\frak q},} and denote by $\mathcal X_{\frak q^\lambda}^{\gc,\sigma},\mathcal X_{\frak q^\lambda}^{\agc,\sigma}$ the induced vector fields on $W^\sigma$ and $W^\sigma/S^1$, respectively.

\subsection{Sobolev completions}
Fix an integer $k\ge 5$. Taking $L^2_k$ Sobolev completions of the configuration spaces and gauge group, we obtain spaces \eq{
W_k\subset\mathcal C_k(Y),\quad W_k^\sigma\subset\mathcal C_k^\sigma(Y),\quad \mathcal G_k(Y),\quad \mathcal B_k^\sigma(Y)=\mathcal C_k^\sigma(Y)/\mathcal G_{k+1}(Y).}
Let $\mathcal T_j^{\gc}$ be the $L^2_j$-completion of the tangent bundle of $W_k$, namely the trivial bundle over $W_k$ with fiber $W_j$. The bundle $\mathcal T_j^{\gc}$ extends to a bundle $\mathcal T_j^{\gc,\sigma}$ over  $W^\sigma$, and the anticircular global Coulomb slices complete to a subbundle $\mathcal K^{\agc,\sigma}_j\subset \mathcal T_j^{\gc,\sigma}$. The vector field $\mathcal X_{\frak q}^{\gc}$ is then a smooth section $W_k\to\mathcal T_{k-1}^{\gc}$.

\subsection{Four-dimensional configurations}\label{subsection:4dconf}
In the remaining subsections, we consider the cylinder $Z=I\times Y$, where $Y$ is a homology 3-sphere and $I\subset \R$ is an interval (possibly all of $\R$). We also fix a $k>0$, a $\delta>0$, and a smooth function $g:\R\to\R$ with $g(t)=e^{\delta |t|}$ for $|t|\gg 0$, giving rise to weighted Sobolev norms on $Z$ when $I=\R$.

The basic configuration space for monopole Floer homology is \eq{
    \mathcal C(Z)=\set{(a,\phi)\mid a\in\Omega^1(Z;i\R),\phi\in\Gamma(\mathbb S^+)}.}
We define \eq{
    \mathcal C^\tau(Z)\subset \Omega^1(Z;i\R)\times C^\infty (I)\times C^\infty(Z;\mathbb S^+)}
to be the space of triples $(a,s,\phi)$ with $s(t)\ge 0$, $\|\phi(t)\|_{L^2(Y)}=1$ for all $t\in I$.
$\mathcal C^\tau(Z)$ is a subset of the space $\tilde {\mc C}^\tau(Z)$, defined by removing the condition $s(t)\ge 0$ in the definition of $\mathcal C^\tau(Z)$. Then $\tilde{\mathcal C}^\tau(Z)$ has an $L^2_{k,\delta}$ completion to a Hilbert manifold $\tilde {\mathcal C}^\tau_{k,\delta}(Z)$, which contains the $L^2_{k,\delta}$ completion $\mathcal C^\tau_k(Z)$ as a closed subset. The completed tangent bundle to $\tilde{\mc C}^\tau_{k,\delta}(Z)$ is the $L^2_j$-completion of the bundle $\mc T^\tau(Z)$ with fibers \eq{
    \mc T^\tau(Z)_{(a,s,\phi)}&=\set{(b,r,\psi)\mid \text{Re}\inner{\phi(t),\psi(t)}_{L^2(Y)}=0,\forall t}\\
    &\subset C^\infty(Z;iT^\ast Z)\oplus C^\infty([t_1,t_2])\oplus C^\infty(Z;\bb S^+).}

The trivial $\Gamma(Z; i\Lambda^2_+T^\ast Z\oplus\bb S^-)$-bundle $\mathcal V(Z)$ over $\mathcal C(Z)$ produces a bundle $\mc V^\tau(Z)$ over $\mc C^\tau(Z)$ with fiber \eq{
    \mc V^\tau(Z)_{(a,s,\phi)}&=\set{(b,r,\psi)\mid \text{Re}\inner{\phi(t),\psi(t)}_{L^2(Y)}=0,\forall t}\\
    &\subset C^\infty(Z;i\Lambda^2_+T^\ast Z)\oplus C^\infty(\R)\oplus C^\infty(Z;\mathbb S^-).}
The bundle $\mathcal V^\tau(Z)$ extends naturally over $\tilde{\mc C}^\tau(Z)$, and the four-dimensional Seiberg-Witten equations give a section $\mc F^\tau_{\frak q}$ of $\mc V^\tau_{k-1,\delta}(Z)$ over $\tilde{\mc C}_{k,\delta}^\tau(Z)$.
\subsection{Gauge-fixing on cylinders}\label{Gauge-fixing on cylinders}
\begin{definition}
We say $(a(t)+\alpha(t)\,dt,\phi(t))\in\Omega^1(Z;i\R)\oplus \Gamma(\bb S^+)$ is in \textit{pseudo-temporal gauge} if $\alpha(t)$ is constant on $Y$ for each fixed $t$.
\end{definition}
We define $\mc C^\gc(Z)\subset \mc C(Z)$ to consist of configurations $(a(t)+\alpha(t)\,dt,\phi(t))$ in pseudo-temporal gauge for which $(a(t),\phi(t))$ is in slicewise Coulomb gauge, i.e. \eq{
    d(\alpha(t))=0,\; d^\ast(a(t))=0,\;\forall t.}
Inside $\mc C^\gc(Z)$ is the subset $W(Z)$ of configurations in temporal gauge, i.e. $\alpha=0$. (For each of these spaces, we can also consider a $\sigma$ or $\tau$ blowup and the embedding of the $\tau$ blowup into a Hilbert manifold, just as we did with $\mc C(Z)$.) In this setting, the Seiberg-Witten equations can be thought of as a section \eq{
    \mc F^{\gc}_{\frak q}:\mc C^{\gc}(Z)\to\mc V^{\gc}(Z),} where $\mc V^{\gc}(Z)$ is the trivial $W(Z)$-bundle over $\mc C^{\gc}(Z)$.
    
\subsection{Path spaces}
Given points $x,y\in W^\sigma$ and a smooth path $\gamma_0:\R\to W^\sigma$ such that $\gamma_0(t)=x$ for $t\ll 0$, $\gamma_0(t)=y$ for $t\gg 0$, write $Z=\R\times Y$ and define \eq{
    \mc C^{\gc,\tau}_k(x,y)=\set{\gamma\in \mc C^{\gc,\tau}_{k,\loc}(\R\times Y)\mid \gamma-\gamma_0\in L^2_k(Z; iT^\ast Z)\times L^2_k(\R;\R)\oplus L^2_k(Z;\bb S^+)}.}
This has an action of the gauge group \eq{
    \mc G^{\gc}_{k+1}(Z)=\set{u:\R\to S^1\mid 1-u\in L^2_{k+1}(\R;\C)},} and we denote the quotient by $\mc B^{\gc,\tau}_k([x],[y])$. 
    Furthermore, $\mc C^{\gc,\tau}_k(x,y)$ embeds in a Hilbert manifold $\tilde{\mc C}^{\gc,\tau}_k(x,y)$ 
    (obtained by dropping the condition $s(t)\ge 0$), and we write \eq{
    \tilde{\mc B}^{\gc,\tau}_k([x],[y])=\tilde{\mc C}^{\gc,\tau}_k(x,y)/\mc G^{\gc}_{k+1}(Z).}
    We also set
\eq{
    W^\tau_k(x,y)=\set{\gamma\in \mc C^{\gc,\tau}_k(x,y)\mid \gamma\text{ is in temporal gauge}}.}

\section{Gradings}\label{sec:gradings}
In this short section, we define the relative gradings between two stationary points of the Seiberg-Witten flow, and show that these gradings can be expressed as the Fredholm index of a certain operator on the global Coulomb slice alone.

\subsection{The Morse-Bott index}\label{subsection:mb-index}
This subsection follows section 3 in Chapter 2 of \cite{lin}.

Choose two contractible open sets $[\frak U_-],[\frak U_+]$ inside critical submanifolds $[\frak C_-],[\frak C_+]\subset \mc C^\sigma_k(Y)/\mc{G}$ admitting smooth lifts $\frak U_-,\frak U_+\subset \mc C^\sigma_k(Y)$ contained in local Coulomb slices at $x_\pm\in \frak U_{\pm}$. There is then a Banach manifold $\tilde{\mc C}^\tau_{k,\delta}(\frak U_-,\frak U_+)$ consisting of paths from $\frak U_-$ to $\frak U_+$ whose difference from some reference path from $\frak U_-$ to $\frak U_+$ is in $L^2_{k,\delta}$. Next, we define $T^\tau_{j,\delta}$ to be the $L^2_{j,\delta}$ completion of the tangent bundle of $\tilde{\mc C}^\tau_{k,\delta}(\frak U_-,\frak U_+)$ (see \cite{lin} for details) and set $T^\tau_{j,\delta,0}\subset T^\tau_{j,\delta}$ to be the finite-codimension sub-bundle consisting of vector fields vanishing at the endpoints. We set \eq{
    \mathbf d^\tau:\text{Lie}(\mc G_{j+1,\delta})\times \tilde{\mc C}^\tau_{k,\delta}(\frak U_-,\frak U_+)\to \mc T^\tau_{j,\delta},\\
    \mathbf d^\tau(\xi,\gamma)=(-d\xi,0,\xi\phi_0),} where $\gamma=(a_0,s_0,\phi_0)$, to be the linearization of the gauge group action.
    We also define an operator \eq{
    \mathbf d^{\tau,\dagger}:\mc T^\tau_{j,\delta}\to L^2_{j-1,\delta}(Z;i\R)} given by \eq{
    (a,s,\phi)\mapsto -d^\ast a-2\sigma(t)c+is_0^2\text{Re}\inner{i\phi_0,\phi}+i|\phi_0|^2\text{Re}(\mu_Y(\inner{i\phi_0,\phi})),}
    where $a=c\,dt + b$, $\mu_Y$ is average value over $Y$, and $\sigma(t)=g'/g$, where $g$ is the function defining the weighted Sobolev norm.
    \begin{remark}
    This definition of $\mathbf d^{\tau,\dagger}$ differs from the definition given in \cite{lm}, equation (5.8.65), by the addition of the term $-2\sigma(t)c$. This term has to be included in the Morse-Bott setting because we are working with weighted spaces; the map $a\mapsto -d^\ast a-2\sigma(t)c$ is the adjoint of $-d$ with respect to the weighted norms.
    \end{remark}
    Now we define, for $\gamma\in\tilde{\mc C}^\tau_k(\frak U_-,\frak U_+)$, \eq{
    \mc K^\tau_{j,\delta,\gamma}&=\ker \mathbf d_\gamma^{\tau,\dagger}\subset \mc T^\tau_{j,\delta,\gamma}\\
    \mc J^\tau_{j,\delta,\gamma}&=\im\mathbf d_\gamma^{\tau}\subset \mc T^\tau_{j,\delta,\gamma}.}
    
To differentiate the section $\mc F^\tau_{\frak q}$ of the bundle $\mc V^\tau_{k-1,\delta}\to\mc C^\tau_{k,\delta}(\frak U_-,\frak U_+)$, we first define a projection \eq{
    \Pi_\gamma^\tau:L^2_{j,\delta}(Z;i\frak{su}(\bb S^+))\oplus L^2_{j,\delta}(\R;\R)\oplus L^2_{j,\delta}(Z;\bb S^-)\to \mc V^\tau_{j,\delta,\gamma}} by applying $L^2$ projection on each slice $\{t\}\times Y$. Then we set $\mc D\mc F^\tau_{\frak q}$ to be the derivative of $\mc F^\tau_{\frak q}$ in the ambient Hilbert space, followed by $\Pi_\gamma^\tau$.

Given a path $\gamma_0:\R\to\mc C^\sigma(Y)$, define an operator \eq{
    Q_{\gamma_0}=\mc D_{\gamma_0} \mc F^\tau_{\frak q}\oplus \mathbf d^{\tau,\dagger}_{\gamma_0}:\mc T^\tau_{j,\delta,\gamma_0,0}\to \mc V^\tau_{j-1,\delta,\gamma_0}\oplus L^2_{j-1,\delta}(Z;i\R).}
\begin{prop}[Proposition 3.4 in Chapter 2 of \cite{lin}]
Suppose we are given two critical submanifolds $\frak U_-, \frak U_+$ in Coulomb slices. Then for each $\gamma_0\in \mc C^\tau_{k,\delta}(\frak U_-,\frak U_+)$, the linear operator \eq{
     Q_{\gamma_0}:\mc T^\tau_{j,\delta,\gamma_0,0}\to \mc V^\tau_{j-1,\delta,\gamma_0}\oplus L^2_{j-1,\delta}(Z;i\R)} above is Fredholm for every $1\le j\le k$. The index of $Q_{\gamma_0}$ is independent of $j$ and $\delta>0$ sufficiently small.
\label{basic-fredholm-weighted}
\end{prop}
\begin{cor}
The restriction of the bundle map \eq{
    \mc D\mc F^\tau_{\frak q}: \mc K^\tau_{j,\delta,\gamma,0}\to \mc V^\tau_{j-1,\delta,\gamma}}
is Fredholm and has the same index as $Q_\gamma$.
\end{cor}

\begin{definition}
Given critical points $x_-,x_+\in\mc C^\sigma_k(Y)$ with images in $[\frak C_-]$ and $[\frak C_+]$, we define \eq{
    \gr(x_-,x_+)=\text{ind}(Q_\gamma)+\dim [\frak C_+]} for any path $\gamma$ from $x_-$ to $x_+$. Given $[x_-],[x_+]\in \mc B^\sigma_k(Y)$ and a relative homotopy class $z\in\pi_1(\mc B^\sigma_k(Y),[\frak C_-],[\frak C_+])$ connecting them, we also define \eq{
        \gr_z([x_-],[x_+])&=\gr(x_-,x_+)} for any pair of lifts such that a path connecting them represents the given homotopy class in the quotient.
\end{definition}

There is also an absolute rational grading defined for a critical point $[x]$ with image in $[\frak C_+]$, given by \eqq{
    \gr^{\Q}([x])=-\gr_z([\mathfrak C_0],X,[\mathfrak C])+\frac{c_1(\mathfrak t)^2-\sigma(X)}4-\iota(X)\in\Q,
    \label{eqn:abs-rat-grad}} where: 
    \begin{itemize}
    \item $X$ is a cobordism from $S^3$ to $Y$;
    \item $[\mathfrak C_0]$ is the reducible critical manifold on $S^3$ which corresponds to the smallest positive eigenvalue of the Dirac operator;
    \item $z$ is any homotopy class connecting critical points in $[\mathfrak C_0]$ and $[\mathfrak C]$;
    \item $\gr_z$ is defined as in \cite{lin}, Section 3.6;
    \item $\mathfrak t$ is a spin$^c$-structure on $X$ extending $\mathfrak s$;
    \item $\iota(X)=\frac12(\chi(X)+\sigma(X)-b_1(Y))$.
    \end{itemize}
\subsection{Restriction to the global Coulomb slice}
We now wish to repeat the constructions of the previous section in the global Coulomb slice, following \S5.9 in \cite{lm}. Recall from Section \ref{Gauge-fixing on cylinders} that $\tilde{\mc C}^{\gc,\tau}$ denotes the $\tau$-model of the blowup of trajectories in pseudo-temporal gauge, and $W^\tau$ denotes those additionally in temporal gauge. There is a section \eq{
    \mc F^{\gc,\tau}_{\frak q}:\tilde {\mc C}^{\gc,\tau}_{k,\delta}(x,y)\to \mc V^{\gc,\tau}_{k-1,\delta}(Z).}

We define \eq{
    \mc K^{\gc,\tau}_{j,\delta,\gamma}=\ker\mathbf d^{\gc,\tau,\tilde{\dagger}}_\gamma\subset \mc T^{\gc,\tau}_{j,\delta,\gamma},} where if we write $b=\beta\,dt+b(t)$,\eq{
    \mathbf d^{\gc,\tau,\tilde{\dagger}}_\gamma(b,r,\psi)=\frac{d\beta}{dt}+\sigma(t)\beta(t)+i\inner{i\phi,\psi}_{\tilde g}.}
(This also differs from the definition of $\mathbf d^{\gc,\tau,\tilde{\dagger}}_\gamma$ in equation (5.8.67) in \cite{lm}, because we again want to compute the adjoint with respect to the weighted spaces.) Then we have the following two propositions, which are proved in the same way as Proposition 5.40 and Proposition 5.42 in \cite{lm}, substituting the weighted Sobolev spaces throughout.
\begin{prop}
Let $x,y$ be normally hyperbolic stationary points of $\mc X^\sigma_{\frak q}$. Pick a path $\gamma\in W^\tau_{k,\delta}(x,y)$. Then for $j\le k$, the operator \eq{
    \left(\mc D_\gamma^\tau \mc F^{\gc,\tau}_{\frak q}\right)|_{\mc K^{\gc,\tau}_{j,\delta,\gamma}}:\mc K^{\gc,\tau}_{j,\delta,\gamma}\to \mc V^{\gc,\tau}_{j-1,\delta,\gamma}}
is Fredholm. Moreover, the Fredholm index is the same as that of the operator \eq{
    \left(\mc D_\gamma^\tau \mc F^{\tau}_{\frak q}\right)|_{\mc K^{\tau}_{j,\delta,\gamma}}:\mc K^{\tau}_{j,\delta,\gamma}\to \mc V^{\tau}_{j-1,\delta,\gamma}.}

\label{coul-gradings}
\end{prop}
\begin{prop}
Consider a path $\gamma\in C^\tau_{k,\delta}(x,y)$ in temporal gauge. Write $x^\flat=\Pi^{\gc,\sigma}(x)$, $y^\flat=\Pi^{\gc,\sigma}(y)$, $\gamma^\flat=\Pi^{\gc,\tau}(\gamma)$. Then:
\begin{enumerate}[(a)]
\item The operators \eq{
   \left(\mc D_\gamma^\tau \mc F^{\tau}_{\frak q}\right)|_{\mc K^{\tau}_{j,\delta,\gamma}}:\mc K^{\tau}_{j,\delta,\gamma}\to \mc V^{\tau}_{j-1,\delta,\gamma}\text{ and } \left(\mc D_{\gamma^\flat}^\tau \mc F^{\gc,\tau}_{\frak q}\right)|_{\mc K^{\gc,\tau}_{j,\delta,\gamma^\flat}}:\mc K^{\gc,\tau}_{j,\delta,\gamma^\flat}\to \mc V^{\gc,\tau}_{j-1,\delta,\gamma^\flat}} have the same Fredholm index.
\item Suppose that $\gamma$ is a trajectory of $\mc X_{\frak q}^\sigma$, so that $[\gamma^\flat]\in\mc B^{\gc,\tau}_k([x^\flat],[y^\flat])$ is a trajectory of $\mc X_{\frak q}^{\agc,\sigma}$. If $ \left(\mc D_\gamma^\tau \mc F^{\tau}_{\frak q}\right)|_{\mc K^{\tau}_{j,\delta,\gamma}}$ is surjective, then so is $\left(\mc D_{\gamma^\flat}^\tau \mc F^{\gc,\tau}_{\frak q}\right)|_{\mc K^{\gc,\tau}_{j,\delta,\gamma^\flat}}$.
\end{enumerate}

\label{upgrade-coul-gradings}
\end{prop}

Thus, we have expressed the grading between stationary points entirely in terms of configurations in global Coulomb gauge.

\section{Stationary points} \label{stationary-points}

\subsection{Reducible stationary points}

For $a\in W_k$, define $D_{\mathfrak q,a}:L^2_k(\mathbb S)\to L^2_{k-1}(\mathbb S)$ by \eq{
D_{\mathfrak q,a}(\phi)=\mathcal D_{(a,0)} \left(\mathcal X_{\mathfrak q}^{\gc}\right)^1(0,\phi),} where the superscript $1$ represents the spinor component. The operator $D_{\mathfrak q,a}$ is a zero-order perturbation of the Dirac operator, which at a reducible critical point $(a,0)$ is the spinorial part of the Hessian of the perturbed Chern-Simons-Dirac functional. This establishes the following proposition:
\begin{prop}
The operator $D_{\mathfrak q,a}$ is Fredholm. If $(a,0)$ is a reducible stationary point of $\mathcal X_{\mathfrak q}^{\gc}$, then $D_{\mathfrak q,a}$ is self-adjoint in the $L^2$-metric.
\label{basic-dirac-props}
\end{prop}
The importance of the operator $D_{\mathfrak q,a}$ is as follows: 
\begin{prop}
Fix a reducible stationary point $(a,0)$ of the vector field $\mathcal X_{\mathfrak q}^{\gc}$ on $W_k$. Then $(a,0,\phi)\in W^\sigma$ is a stationary point of $\mathcal X_{\mathfrak q}^{\gc,\sigma}$ if and only if $\phi$ is an eigenvector of $D_{\mathfrak q,a}$.
\begin{itemize}
\item If $a\ne 0$, the spectrum of $D_{\mathfrak q,a}$ is simple.
\item If $a=0$, $D_{\mathfrak q,a}$ is $j$-invariant, and its eigenspaces have complex dimension 2.
\end{itemize}
\label{reducible-descr-dirac}
\end{prop}
\begin{proof}
According to equation (5.4.26) in \cite{lm}, $(a,0,\phi)$ is a stationary point of $\mathcal X_{\mathfrak q}^{\gc,\sigma}$ if and only if $\phi$ is an eigenvector of the operator $D_{\mathfrak q,a}$. Since global Coulomb projection preserves $\set0\times L^2_k(Y;S)$ and $\mathcal X_{\mathfrak q}^{\gc,\sigma}$ is $\jmath$-invariant, $D_{\mathfrak q,0}$ is $j$-invariant. The fact that the eigenspaces have complex dimension 1 if $a\ne 0$ and 2 if $a=0$ follows from the fact that $\mathfrak q$ is a regular perturbation, so the stationary points are $\pin(2)$-non-degenerate (Definition 5.2.3 in \cite{lin}).
\end{proof}
Thus, we have the following description of the stationary submanifolds in the quotient $W^\sigma_k/S^1$:
\begin{cor}
The connected stationary submanifolds of the vector field $\mathcal X_{\mathfrak q}^{\agc,\sigma}$ on $W_k^\sigma/S^1$ are
\begin{enumerate}[(1)]
\item isolated points;
\item copies of $\C P^1$ in the reducible locus which are labeled by eigenvalues of $D_{\mathfrak q,0}$.
\end{enumerate}
\label{stationary-types}
\end{cor}
It will also be necessary to work with the approximate operator $D_{\mathfrak q^\lambda,a}:L^2_k(\mathbb S)\to L^2_{k-1}(\bb S)$ defined by the analogous formula \eq{
    D_{\mathfrak q^\lambda,a}(\phi)=\mathcal D_{(a,0)} \left(\mathcal X_{\mathfrak q^\lambda}^{\gc}\right)^1(0,\phi).} 
The operator $D_{\mathfrak q^\lambda,a}$ is no longer part of the Hessian of any functional on $\mathcal C(Y)$, and therefore may no longer be self-adjoint. Nevertheless, we still have the following weaker version of Proposition \ref{basic-dirac-props}:
\begin{prop}
The operator $D_{\mathfrak q^\lambda,a}$ has real eigenvalues.
\label{approximate-dirac-real}
\end{prop}
\begin{proof}
Let $V^\lambda\subset V=L^2_k(\bb S)$ denote the span of the eigenvectors of $D$ with eigenvalues in $(-\lambda,\lambda)$.
As in equation (6.3.8) in \cite{lm}, we may write $D_{q^\lambda,0}=D+(p^\lambda)^1L$, where $L:V\to V$ is self-adjoint and $(p^\lambda)^1$ denotes the spinor component of $p^\lambda$. We may write $(p^\lambda)^1$ as a convex combination \eq{
    (p^\lambda)^1=\sum_{i=1}^m \beta_i \,p^{\lambda_i^\bullet}_V,} where $p^{\lambda}_V:V\to V^\lambda$ is $L^2$-orthogonal projection, $\sum_{i=1}^m \beta_i=1$, and $\beta_m\ne 0$. In an eigenbasis for $D$, the matrix $P$ of $(p^\lambda)^1$ is diagonal, and its unique diagonal entries are of the form $\sum_{i=j}^m \beta_i$ (for some $1\le j\le m$); these entries are all nonzero because $\beta_m\ne 0$. Set $A=P^{-1/2}$ (another diagonal matrix); then since $AP=A^{-1}$, the matrix $APLA^{-1}=A^{-1}LA^{-1}$ is self-adjoint: \eq{
        (A^{-1}LA^{-1})^\ast=(A^\ast)^{-1}L^\ast(A^\ast)^{-1}=A^{-1}LA^{-1}.} Since $A$ commutes with $D$, we thus have that \eq{
        D_{\mathfrak q^\lambda,a}|_{V^\lambda}=A^{-1}(D+APLA^{-1})A} is conjugate to a self-adjoint operator, and hence has real eigenvalues. 
        It immediately follows that $D_{\mathfrak q^\lambda,a}$ has real eigenvalues, since it decomposes as an upper-triangular block matrix \eq{
        \begin{pmatrix} D_{\mathfrak q^\lambda,a}|_{V^\lambda}&\ast\\
        0&D \end{pmatrix}} with respect to the decomposition $V=V^\lambda\oplus (V^\lambda)^\perp$.
\end{proof}
\subsection{Identifying stationary points}\label{Identifying stationary points}
In this section, which is inspired by Chapter 7 of \cite{lm}, we fix:
\begin{itemize}
\item a bound $R>0$ such that all stationary points of $\mathcal X_{\mathfrak q}^{\gc}$ are contained in $B(2R)\subset W_k$
\item a grading range $[-N,N]$, a closed neighborhood $\mathcal N$ and an open neighborhood $\mathcal U\subset\mathcal N$ of the set of stationary points of $\mathcal X_{\mathfrak q}^{\gc,\sigma}$ in this grading range. We assume that the projection of $\mathcal N$ to the blow-down is contained in $B(2R)$, and that $N$ is chosen large enough to contain each reducible stationary point $(a,0,\phi)$ of $\mathcal X_{\mathfrak q}^{\gc,\sigma}$ where $\phi$ is an eigenvector of $D_{\mathfrak q,a}$ with smallest positive eigenvalue.
\item an \textit{eigenvector list}: for each eigenvalue $\mu$ of $D_{\mathfrak q,0}$, we choose an eigenvector $\phi_\mu$ of $D_{\mathfrak q,0}$ with eigenvalue $\mu$.
\end{itemize}
Let $\mathfrak C_{\mathcal N}$, $\mathfrak C_{\mathcal N}^\lambda$ denote the sets of stationary points of $\mathcal X_{\mathfrak q}^{\agc,\sigma}$ and $\mathcal X_{\mathfrak q^\lambda}^{\agc,\sigma}$ (respectively) which live in $\mathcal N/S^1$. The purpose of this section is to prove the following proposition:
\begin{prop}
For $\lambda\gg 0$, there is a one-to-one correspondence \eq{
    \Xi_\lambda:\mathfrak C_{\mathcal N}^\lambda\to\mathfrak C_{\mathcal N}.} This correspondence preserves the type of stationary point (irreducible, stable, unstable, isolated).
\label{point-corr-prop}
\end{prop}
This proposition is proved as Corollary 7.8 in \cite{lm} in the non-equivariant setting. The proof carries over to the equivariant setting for \textit{isolated} stationary points, where the Hessian is nondegenerate. Therefore, it only remains to establish a correspondence between the non-isolated critical points of $\mathcal X_{\mathfrak q}^{\agc,\sigma}$ and $\mathcal X_{\mathfrak q^\lambda}^{\agc,\sigma}$.

We begin by restating a useful lemma from \cite{lm} which still applies in the current $\pin(2)$-equivariant setting:
\begin{lemma}[Lemma 7.2 in \cite{lm}]
Fix $\epsilon>0$. There exists $b\gg 0$ such that for all $\lambda>b$ the following is true. If $x\in\mathcal N\subset W^\sigma$ is a zero of $\mc X_{\mathfrak q^\lambda}^{\gc,\sigma}$, then there exists $x'\in\mathcal N$ such that $\mathcal X_{\mathfrak q}^{\gc,\sigma}(x')=0$ and $x,x'$ have $L^2_k$-distance at most $\epsilon$ in $L^2_k(Y;iT^\ast Y)\oplus \R\oplus L^2_k(Y;\mathbb S)$.
\label{lemma-7.2}
\end{lemma}
A consequence of Lemma \ref{lemma-7.2} is the following:
\begin{prop}[Corollary 7.5 in \cite{lm}]
For $\lambda\gg 0$, all the stationary points of $\mathcal X_{\mathfrak q^\lambda}^{\gc,\sigma}$ in $\mathcal N$ live inside the finite-dimensional blow-up $(W^\lambda)^\sigma$.
\label{cor-7.5}
\end{prop}

To prepare the proof of Proposition \ref{point-corr-prop}, index the eigenvalues of the linear operator $l$ by the condition $|\lambda_n|\le |\lambda_{n+1}|$. Pick a strictly decreasing homeomorphism $f:(0,\infty]\to[0,1)$ with the following properties: 
\begin{itemize}
\item The restriction of $f$ to $(0,\infty)$ is a diffeomorphism onto $(0,1)$;
\item $\lim_{n\to\infty} |\lambda_n|^2 f(|\lambda_{n+1}|)=\infty$.
\end{itemize}
\begin{lemma}[Lemma 7.6 in \cite{lm}]
The map \eq{
    h:W_k\times (-1,1)\to W_{k-1},\quad h(x,r)=x-p^{f^{-1}(|r|)}(x)}
is continuously differentiable, with $\mathcal Dh_{(x,0)}(0,1)=0$ for all $x$.
\label{differentiability-homeomorphism}
\end{lemma}

\begin{prop}
Fix an eigenvalue $\mu_0$ of $D_{\mathfrak q,0}:L^2_k(\bb S)\to L^2_{k-1}(\bb S)$ and choose a corresponding eigenvector $\phi_0$. Then for sufficiently large $\lambda_0$, there are continuous functions \eq{
\phi:(\lambda_0,\infty]\to L^2_k(\bb S),\qquad \mu:(\lambda_0,\infty]\to \bb R} satisfying $\phi(\infty)=\phi_0$, $\|\phi(\lambda)\|_{L^2}=1$, $\mu(\infty)=\mu_0$, and $D_{\mathfrak q^{\lambda},0}(\phi(\lambda))=\mu(\lambda)\phi(\lambda)$ for $\lambda\in (\lambda_0,\infty]$. 

The functions $\phi(\lambda)$ and $\mu(\lambda)$ are unique in the following sense: there is a neighborhood $U_1\subset \R$ of $\mu_0$ and a neighborhood $U_2\subset L^2_k(\mathbb S)$ of $\phi_0$ such that for all $\lambda>\lambda_0$, $\phi(\lambda)$ is the unique unit eigenvector of $D_{\mathfrak q^\lambda,0}$ in $U_2$ (up to multiplication by a unit quaternion) with eigenvalue in $U_1$.
\label{imp-fn-stat}
\end{prop}
\begin{proof}
Let $\iota: L^2_k(\bb S)\to L^2_{k-1}(\bb S)$ be the inclusion, and let $\inner{\cdot,\cdot}_{\bb H}$ be a quaternionic $L^2$-inner product on $\Gamma(\bb S)$. Write $D_r=D_{\mathfrak q^{f^{-1}(|r|)},0}$, and consider the map \eq{
    S:(-1,1)\times L^2_k(\bb S)\times \bb H\to L^2_{k-1}(\bb S)\oplus \bb H\\
    S(r,\phi,\mu)=\left((D_r-\mu \iota)\phi,\inner{\phi,\phi_0}_{\bb H}\right).}
    $S$ is obviously smooth in $\phi$ and $\mu$, and is smooth in $r$ away from $r=0$, since the $p^\lambda$ are smoothed projections. To show differentiability at $r=0$, we compute \eq{
    S(r,\phi,\mu)-S(0,\phi,\mu)&=(D_r\phi-D_0\phi,0).} But \eq{
    D_r\phi-D_0\phi&=\left(D\phi+p^{f^{-1}(|r|)}\left[\mathcal D_{(0,0)}(c_{\mathfrak q})^1(0,\phi)\right]\right)-\left(D\phi+\mathcal D_{(0,0)}(c_{\mathfrak q})^1(0,\phi)\right)\\
    &=p^{f^{-1}(|r|)}\left[\mathcal D_{(0,0)}(c_{\mathfrak q})^1(0,\phi)\right]-\mathcal D_{(0,0)}(c_{\mathfrak q})^1(0,\phi),} so by Lemma \ref{differentiability-homeomorphism}, we have $\mathcal DS_{(0,\phi,\mu)}(1,0,0)=0$.
    
    Write $S_r=S(r,\cdot,\cdot)$. Then we compute \eq{
    \mathcal D(S_0)_{(\phi_0,\mu_0)}{(\psi,\nu)}&=\left((D_0-\mu_0\iota)\psi-\nu\phi_0,\inner{\psi,\phi_0}_{\bb H}\right).} If $\mathcal D(S_0)_{(\phi_0,\mu_0)}(\psi,\nu)=(0,0)$ for some $(\psi,\nu)$, then $(D_0-\mu_0\iota)\psi$ is a multiple of $\phi_0$ and hence lies in the $\mu_0$-eigenspace of $D_0$. Since $D_0$ diagonalizes over $\bb H$ with simple eigenvalues, $(D_0-\mu_0\iota)\psi$ lies in the $\mu_0$-eigenspace of $D_0$ only when $\psi$ itself lies in the $\mu_0$-eigenspace of $D_0$. Since the $\mu_0$-eigenspace is exactly the $\bb H$-span of $\phi_0$, the second equation $\inner{\psi,\phi_0}_{\bb H}=0$ implies $\psi=0$.  Hence $\nu=0$ as well, and $\mathcal D(S_0)_{(\phi_0,\mu_0)}$ is injective. Since $\mathcal D(S_0)_{(\phi_0,\mu_0)}$ is a compact perturbation of the self-adjoint Fredholm operator $D_0\oplus 0:L_k^2(\bb S)\times \bb H\to L_{k-1}^2(\bb S)\times \bb H$, its index must be 0 and hence it is also surjective. By the implicit function theorem, we can therefore uniquely solve the equation \eq{
    S(r,\phi,\mu)=(0,1)} for $\phi$ and $\mu$ in terms of $r$ when $r$ is small. Finally, the image of $\mu$ actually lands in $\R\subset \bb H$ by Proposition \ref{approximate-dirac-real}, and we can replace $\phi(r)$ by $\phi(r)/\|\phi(r)\|_{L^2}$ to make it a unit eigenvector.

\end{proof}
\begin{remark}
Among all eigenvectors of $D_{\mathfrak q^\lambda,0}$ with eigenvalue $\mu(\lambda)$, the eigenvector $\phi(\lambda)$ constructed above is the unique eigenvector closest to $\phi_0$. This follows from the equation $\inner{\phi(\lambda),\phi_0}_{\bb H}>0$, and the fact that if $\phi,\psi$ are any two unit vectors in $\bb H^\infty$, the closest point to $\phi$ on $S(\bb H)\psi$ is $\psi$ if and only if $\inner{\phi,\psi}_{\bb H}>0$.
\end{remark}
\begin{cor}
For $\lambda\gg0$, any reducible stationary point $x=(0,0,\phi)$ of $\mc X_{\mathfrak q^\lambda}^{\agc,\sigma}$ contained in $\mc N$ corresponds to a nonzero eigenvalue of $D_{\mathfrak q^{\lambda},0}$ with complex multiplicity exactly 2.
\label{mult-2-in-N}
\end{cor}
\begin{proof}
It follows from Lemma \ref{lemma-7.2} that there is a reducible stationary point $x_0=(0,0,\phi_0)$ of $\mc X_{\mathfrak q}^{\agc,\sigma}$ arbitrarily close to $x$, and therefore the spinorial energy $\mu_0$ of $x_0$ is arbitrarily close to the $\lambda$-spinorial energy $\mu$ of $x$. We may thus assume that $\mu$ lies within the neighborhood $U_1$ of $\mu_0$ constructed in Proposition \ref{imp-fn-stat}. By the uniqueness statement in Proposition \ref{imp-fn-stat}, it follows that $\mu=\mu(\lambda)$ and $\phi$ is a unit quaternion multiple of $\phi(\lambda)$. If there is an eigenvector $\phi'$ of $D_{\mathfrak q^\lambda,0}$ which is $\bb H$-linearly independent of $\phi$, then we obtain a one-parameter family $[\phi+t\phi']$ of eigenvectors of $D_{\mathfrak q^\lambda,0}$, contradicting the uniqueness statement in Proposition \ref{imp-fn-stat}.
\end{proof}

\begin{proof}[Proof of Proposition \ref{point-corr-prop}]
For isolated stationary points, the correspondence $\Xi_\lambda$ is given by Corollary 7.8 in \cite{lm}. 

Fix a connected stationary locus $B^{\mu_0}$ of $\mathcal X_{\mathfrak q}^{\agc,\sigma}$ corresponding to the projectivization of the $\mu_0$-eigenspace of $D_{\mathfrak q,0}$. Given a $\mu_0$-eigenvector $\phi_0$ of $D_{\mathfrak q,0}$, Proposition \ref{imp-fn-stat} produces an eigenvector $\phi(\lambda)$ of $D_{\mathfrak q^\lambda,0}$ with eigenvalue $\mu(\lambda)$ for every $\lambda\gg 0$. There is then a connected stationary locus $B^{\mu(\lambda)}$ of $\mathcal X_{\mathfrak q^\lambda}^{\agc,\sigma}$ corresponding to the projectivization of $\vspan_{\bb H}(\phi(\lambda))$. We obtain a canonical map $B^{\mu_0}\to B^{\mu(\lambda)}$ induced by the unique $\bb H$-linear map $\vspan_{\bb H}(\phi_0)\to\vspan_{\bb H}(\phi(\lambda))$ which sends $\phi_0$ to $\phi(\lambda)$. By the remark after Proposition \ref{point-corr-prop}, this map is independent of the choice of $\phi_0$.

By repeating this process for every eigenvalue of $D_{\mathfrak q,0}$, we obtain the correspondence $\Xi_\lambda^{-1}$. The proof of Corollary \ref{mult-2-in-N} shows that the correspondence is one-to-one.
\end{proof}

\subsection{Equating relative gradings}
In this section, we prove the following:
\begin{prop}
The correspondence $\Xi_\lambda$ from Proposition \ref{point-corr-prop} preserves relative gradings.
\label{prop:xi-pres-grad}
\end{prop}
\begin{proof}
This follows from Proposition \ref{approx-grading} below.
\end{proof}
From equation (6.3.7) in \cite{lm}, we obtain an approximate section \eq{
    \mc F^{\gc,\tau}_{\mathfrak q^\lambda}:\tilde C^{\gc,\tau}(Z)\to \mc V^{\gc,\tau}(Z).} 
Let $[x_\infty],[y_\infty]$ be stationary points of $\mc X^{\agc,\sigma}_{\mathfrak q}$ in $\mc N/S^1$ and let $[x_\lambda],[y_\lambda]$ be the corresponding stationary points of $\mc X^{\agc,\sigma}_{\mathfrak q^\lambda}$ under the correspondence $\Xi_\lambda$ from Proposition \ref{point-corr-prop}.
\begin{prop}
For $1\le j\le k$ and $\lambda\gg0$, the following is true: for each path $[\gamma_\lambda]\in\mc B^{\gc,\tau}_{k,\delta}([x_\lambda],[y_\lambda])$ with representative $\gamma_\lambda\in \mc C^{\gc,\tau}_{k,\delta}(x_\lambda,y_\lambda)$, the operator \eq{
    \left(\mc D_\gamma^\tau \mc F^{\gc,\tau}_{\mathfrak q^\lambda}\right)|_{\mc K^{\gc,\tau}_{j,\delta,\gamma_\lambda}}:\mc K^{\gc,\tau}_{j,\delta,\gamma_\lambda}\to \mc V^{\gc,\tau}_{j-1,\delta,\gamma_\lambda}} is Fredholm, with index equal to $\gr([x_\infty],[y_\infty])$.
\label{approx-grading}
\end{prop}
\begin{proof}
By Proposition 5.40 in \cite{lm}, $\gr([x_\infty],[y_\infty])$ is equal to the index of an operator \eq{
	\widehat{Q}^{\gc}_{\gamma}:\mc T^\tau_{j,\gamma}\to \mc T^\tau_{j-1,\gamma},} where $\mc T^\tau_{k,\gamma}$ is the pullback over $\gamma$ of the $L^2_k$-completion of the bundle $\mc T^\tau(Z)$ defined in \ref{subsection:4dconf}. Similarly, the index of the operator in the statement of the proposition is equal to the index of an operator \eq{
	\widehat{Q}^{\gc}_{\gamma_\lambda,\mathfrak q^\lambda}:\mc T^\tau_{j,\gamma_\lambda}\to\mc T^\tau_{j-1,\gamma_\lambda}.}
Choose a bump function $\beta(t)$ with $\beta(t)=1$ for $|t|\le T+1$, $\beta(t)=0$ for $|t|\ge T+2$, and set \eq{
	\widehat{Q}^{\text{int},\lambda}=(1-\beta(t))\widehat{Q}^{\gc}_{\gamma_\lambda,\mathfrak q^\lambda}+\beta(t)\widehat{Q}^{\gc}_{\gamma}.}
	The proof of Lemma 9.1 in \cite{lm} shows that for each $t$, the difference  $\displaystyle\widehat{Q}^{\gc}_{\gamma_\lambda,\mathfrak q^\lambda}-\widehat{Q}^{\text{int},\lambda}$ is continuous as a map $L^2_j(Y)\to L^2_{j-1}(Y)$; hence, since $\beta$ is compactly supported, it is continuous as a map $L^2_{j,\delta}(\R\times Y)\to L^2_{j-1,\delta}(\R\times Y)$.
		
	On the other hand, the proof of Lemma 9.1 in \cite{lm} also shows that the difference $\displaystyle \widehat{Q}^{\gc}_{\gamma}-\widehat{Q}^{\text{int},\lambda}$ is equal to $(1-\beta(t))A_\lambda$, where $A_\lambda:L^2_n(Y)\to L^2_{n-1}(Y)$ is an operator which converges to zero in norm as $\lambda\to\infty$. We now wish to conclude that $(1-\beta(t))A_\lambda\to 0$ when considered as an operator $L^2_{j,\delta}(\R\times Y)\to L^2_{j-1,\delta}(\R\times Y)$; however, since $\beta$ is compactly supported, it is enough to show $A_\lambda\to 0$ as an operator between these spaces.
	
	 For $\eta\in L^2_{j,\delta}(\R\times Y)$, \eq{
    \|A_\lambda(\eta(t))\|^2_{L^2_{j-1,\delta}(\R\times Y)} &= \sum_{n=0}^{j-1}\int_{\R} \|A_\lambda(\eta^{(n)}(t))\|^2_{L^2_{j-n-1}(Y)}\,g(t)\,dt\\
    &\le \sum_{n=0}^j \int_{\R} C_\lambda\|\eta^{(n)}(t)\|^2_{L^2_{j-n}(Y)}\,g(t)\,dt\\
    &=C_\lambda \|\eta\|^2_{L^2_{j,\delta}(\R\times Y)},} where $C_\lambda\to 0$. This establishes the desired convergence.
\end{proof}

\subsection{Finite-dimensional gradings}
According to Corollary 7.5 in \cite{lm}, the stationary points of $\mc X^{\agc,\sigma}_{\mathfrak q^\lambda}$ within $\mc N/S^1$ are contained in $(W^\lambda)^\sigma$, and therefore have an associated finite-dimensional Morse-Bott index from the flow of $\mc X^{\agc,\sigma}_{\mathfrak q^\lambda}|_{(B(2R)\cap W^\lambda)^\sigma/S^1}$.
\begin{prop}[Proposition 9.4 in \cite{lm}]
Let $[x_\lambda], [y_\lambda]$ be stationary points of $\mc X^{\agc,\sigma}_{\mathfrak q^\lambda}$. Then $\gr([x_\lambda],[y_\lambda])$, where this grading is computed in infinite dimensions, is equal to the difference in Morse-Bott gradings of $[x_\lambda]$, $[y_\lambda]$ as stationary points of $\mc X^{\agc,\sigma}_{\mathfrak q^\lambda}|_{(B(2R)\cap W^\lambda)^\sigma/S^1}$.
\label{grading-inf-fin}
\end{prop}
\begin{proof}
The proof of Proposition 9.4 in \cite{lm} goes through in the current weighted setting.
\end{proof}

\subsection{Hyperbolicity}
In this section, we show that for large enough $\lambda$, the stationary points of $\mc X_{\mathfrak q^\lambda}^{\agc,\sigma}$ in $\mathcal N/S^1$ have hyperbolic normal Hessians. This is already accomplished by Proposition 7.10 in \cite{lm} for isolated stationary points, where the normal Hessian is just the Hessian. Therefore, we only need to prove hyperbolicity for the non-isolated stationary points, as in the following proposition.

\begin{prop}
For $\lambda\gg 0$, the eigenvalues of the normal Hessian $\hess^{\sigma,\nu}_{\mathfrak q^\lambda,x}$ at a reducible stationary point $x=(0,0,\psi_0)\in \mathcal N$ of $D_{\mathfrak q^\lambda,0}$ are nonzero real numbers.
\label{hyperbolic-reducibles}
\end{prop}
\begin{proof}
The proof of Proposition 5.2.4 in \cite{lin} carries over to the approximate case, giving us an expression \eq{
    \hess^{\sigma}_{\mathfrak q^\lambda,x}=\begin{pmatrix}
    \mu_0 & 0 &0\\
    \ast & \ast d_{\mathfrak q^\lambda} & 0\\
    \ast & \ast & D_{\mathfrak q^\lambda,0}-\mu_0 \end{pmatrix}:\R\oplus \ker d^\ast\oplus (\C\psi_0)^\perp\to\R\oplus \ker d^\ast\oplus (\C\psi_0)^\perp.}
    Here $D_{\mathfrak q^\lambda,0}(\psi_0)=\mu_0\psi_0$ and \eq{
    \ast d_{\mathfrak q^\lambda}(b)=\ast db+2(p^\lambda)^0\mathcal D_0 \mathfrak q^0(b,0),} where the superscript $0$ denotes the form component. 

    Letting $B\subset W^\sigma$ be the critical locus, we have $T_xB=\C(\psi_0\cdot j)$, so the normal Hessian is \eq{
    \hess^{\sigma,\nu}_{\mathfrak q^\lambda,x}=\begin{pmatrix}
    \mu_0 & 0 &0\\
    \ast & \ast d_{\mathfrak q^\lambda} & 0\\
    \ast & \ast & D_{\mathfrak q^\lambda,0}-\mu_0 \end{pmatrix}:\R\oplus \ker d^\ast\oplus (\C\psi_0)^\perp/\C(\psi_0\cdot j)\to\R\oplus \ker d^\ast\oplus (\C\psi_0)^\perp/\C(\psi_0\cdot j).}
    
    The first entry $\mu_0$ is nonzero because the perturbation $\mathfrak q$ is non-degenerate. It follows just as in the proof of Proposition \ref{approximate-dirac-real} that the operator $\ast d_{\mathfrak q^\lambda}$ has real eigenvalues. The statement that $\ast d_{\mathfrak q^\lambda}$ has no zero eigenvalue is equivalent to the statement that $0$ is a non-degenerate zero of $\left(\mc X_{\mathfrak q^\lambda}^{\gc,\sigma}\right)^0$, and this is established in Proposition 7.10 in \cite{lm}. Finally, since the eigenspace of $\mu_0$ is exactly $\vspan_\C(\psi_0,\psi_0\cdot j)$ by Corollary \ref{mult-2-in-N}, it follows that the lower-right block $D_{\mathfrak q^\lambda,0}-\mu_0$ has nonzero real eigenvalues as an operator on $(\C\psi_0)^\perp/\C(\psi_0\cdot j)$.

\end{proof}

\subsection{Other stationary points}
In the previous section, we proved that stationary points of $\mathcal X_{\mathfrak q^\lambda}^{\agc,\sigma}$ in $\mathcal N/S^1$ are normally hyperbolic. However, we are also interested in the stationary points in $B(2R)^\sigma$ which do not lie in $\mathcal N/S^1$. For these stationary points, we obtain the following proposition by combining the proof of Proposition 7.11 in \cite{lm} with the proof of Theorem 5.2.6 in \cite{lin}:

\begin{prop}
We can choose the admissible perturbation $\mathfrak q$ such that for any $\lambda\in\set{\lambda_1^\bullet,\lambda_2^\bullet,\dots}$ sufficiently large, the restriction of $\mathcal X_{\mathfrak q^\lambda}^{\gc,\sigma}$ to $(B(2R)\cap W^\lambda)^\sigma$ has only normally hyperbolic stationary points.
\label{all-normally-hyperbolic}
\end{prop}

\subsection{The quasi-gradient condition}
\begin{prop}
We can choose the perturbation $\mathfrak q$ such that for all $\lambda=\lambda_i^{\bullet}$ with $i\gg 0$, the vector field $\mathcal X_{\mathfrak q^\lambda}^{\gc}$ on $W^\lambda\cap B(2R)$ is a Morse-Bott equivariant quasi-gradient, in the sense of Definition \ref{morse-eq-qg}.
\label{prop:quasi-grad-satisfied}
\end{prop}
\begin{proof}
We must check that $\mathcal X_{\mathfrak q^\lambda}^{\gc}$ satisfies the following conditions:
\begin{enumerate}[(1)]
\item all stationary points of $\mathcal X_{\mathfrak q^\lambda}^{\gc}$ in $W^\lambda\cap B(2R)$ are hyperbolic;
\item the operators $\mathcal D_x\mathcal X_{\mathfrak q^\lambda}^{\gc}$ are self-adjoint;
\item there exists a smooth $\pin(2)$-invariant function $\tilde f: W^\lambda\cap B(2R)\to \R$ such that $d\tilde f(\mathcal X_{\mathfrak q^\lambda}^{\gc})\ge 0$ at all $x\in W^\lambda\cap B(2R)$, with equality holding if and only if $x$ is a stationary point of $\mathcal X_{\mathfrak q^\lambda}^{\gc}$.
\end{enumerate}
Condition (1) is Proposition \ref{all-normally-hyperbolic}. Condition (2) follows from the fact that for $\lambda=\lambda_i^\bullet$, the projection $p^\lambda$ is the $L^2$-orthogonal projection onto $W^\lambda$, and hence \eq{
    \mathcal D_x\mathcal X_{\mathfrak q^\lambda}^{\gc}&=l+p^\lambda \mathcal D_x c_{\mathfrak q}=l+p^\lambda \mathcal D_x c_{\mathfrak q}p^\lambda} is self-adjoint.
As for condition (3), Proposition 8.2 in \cite{lm} provides an $S^1$-invariant function $F_\lambda$ satisfying all necessary properties except $\pin(2)$-invariance. Then $\tilde f=F_\lambda + F_\lambda\circ\jmath$ is $\pin(2)$-invariant and still satisfies the remaining properties.
\end{proof}

\section{Interpolating moduli spaces}\label{sec:inter-mod-sp}

In this section, we describe a chain homotopy equivalence between Morse-Bott complexes built from the Seiberg-Witten flow on the one hand and the approximate flow on the other hand. The construction of a chain homotopy equivalence diverges from \cite{lm}, where it was possible to construct an explicit isomorphism between the two complexes.
\subsection{Identifying the complexes}

Choose $N>0$ large enough that
\begin{enumerate}[(1)]
\item all irreducible stationary points of $\xagcsq$ have grading in $[-N,N]$;
\item no boundary-stable stationary point has grading less than $-N$.
\end{enumerate}
Denote by $\mathfrak C_{[-N,N]}\subset W^\sigma/S^1$ the set of stationary points in this grading range, let $\mathfrak D_{[-N,N]}$ be its preimage in $W^\sigma$, and define \eq{
    \mc N&=\set{x\in W^\sigma_k\mid d_{L^2_k}(x,\mathfrak D_{[-N,N]})<2\delta}.}

Let $\scr S^\lambda\subset B(2R)\subset W$ be the isolated invariant set consisting of points on trajectories of $\xgcql$ which remain in $B(2R)$ for all time. By the following lemma, $\mathscr S^\lambda$ is actually contained in the finite-dimensional subspace $W^\lambda$:
\begin{lemma}[Lemma 6.14 in \cite{lm}]
If $\gamma(t)$ is a trajectory of $\xgcql$ contained in $B(2R)$, then $\gamma(t)$ is contained in $W^\lambda$.
\label{lemma:stays-in-Wlambda}
\end{lemma}
Recall also from Proposition \ref{prop:quasi-grad-satisfied} that for $\lambda=\lambda^{\bullet}_i$, the vector field $\xgcsql$ on $W^\lambda\cap B(2R)^\sigma$ is a Morse-Bott equivariant quasi-gradient. Therefore, for $\lambda=\lambda^{\bullet}_i$, we may follow the construction in \S \ref{subsection:mb-homology-mfld-pin2} to define a finite-dimensional Morse-Bott complex $(\check C_\lambda,\check \partial_\lambda)$ for the vector field $\xagcsql$ on $(W^\lambda\cap B(2R))^\sigma/S^1$, using the isolated invariant set $\scr S^\lambda$.

When $\lambda=\infty$, we define a chain complex $(\check C_\infty,\check \partial_\infty)$ using the flow of $\xagcsq$ in exactly the same way, except that we no longer have finite-dimensional Morse-Bott gradings available. The generators of $(\check C_\infty,\check \partial_\infty)$ are $\delta$-chains $\sigma:\Delta\to B$, where $B\subset B(2R)$ is a stable or interior stationary submanifold of $\xagcsq$; we define the degree of $\sigma$ to be $\dim\Delta+\gr^{\Q}(B)$, where $\gr^{\Q}(B)$ is the absolute rational grading of any point in $B$, as given by (\ref{eqn:abs-rat-grad}). As usual, the differential $\check\partial_\infty$ is defined by the formula (\ref{eqn:bdry-def}). 

Next, we define $\widecheck{\CM}^\lambda(Y,\mathfrak q)$ to be $(\check C_\lambda,\check \partial_\lambda)$ with gradings shifted down by $\dim W^{(-\lambda,0)} + 2n(Y,\mathfrak s, g)$, and define $\widecheck{\HM}^\lambda(Y,\mathfrak q)$ to be the homology of this complex. We write $\widecheck{\CM}(Y,\mathfrak q)=\widecheck{\CM}^\infty(Y,\mathfrak q)$ for the complex $(\check C_\infty,\check \partial_\infty)$, with no grading shift. For $\lambda\le \infty$, we also define $\widecheck{\CS}^\lambda(Y,\mathfrak q)\subset \widecheck{\CM}^\lambda(Y,\mathfrak q)$ to be the subcomplex of $\jmath$-invariant chains, and denote its homology by $\widecheck{\HS}^\lambda(Y,\mathfrak q)$.

The purpose of the grading shifts in the definition of $\widecheck{\CM}^\lambda(Y,\mathfrak q)$ is shown by the following proposition, which is a straightforward extension of \cite{lm}, Proposition 9.9:
\begin{prop}
Let $\sigma:\Delta\to B$ be a generator of $\widecheck{CM}(Y,\mathfrak s)$. For $\lambda\gg 0$, let $\sigma_\lambda=\Xi_{\lambda}^{-1}\circ \sigma$, where $\Xi_\lambda$ is the correspondence defined in Proposition \ref{point-corr-prop}. Then for any $\lambda=\lambda_i^\bullet \gg0$, the grading of $\sigma_\lambda$ in $\widecheck{CM}^\lambda(Y,\mathfrak s)$ agrees with the grading of $\sigma$ in $\widecheck{CM}(Y,\mathfrak s)$.
\end{prop}

Recall from the construction in \S \ref{subsubsection:QV} that for all $\lambda\le\infty$, $\widecheck{\HM}^\lambda$ and $\widecheck{\HS}^\lambda$ are modules over the ring $\mathcal R=\mathbb F[Q,V]/(Q^3)$, where $\deg Q=-1$, $\deg V=-4$, and the action is defined by fiber products with cut-down moduli spaces.

\begin{prop}
For $i\in[-(N-4),N-4]$, there are isomorphisms \eq{
 \widecheck{\HM}^\lambda_i(Y,\mathfrak q) \simeq \widecheck{\HM}_i(Y,\mathfrak q)\\
 \widecheck{\HS}^\lambda_i(Y,\mathfrak q) \simeq \widecheck{\HS}_i(Y,\mathfrak q)}
 respecting the $\mathcal R$-module structure.
 \label{prop:lambda-ultimate}
\end{prop}

The proof of Proposition \ref{prop:lambda-ultimate} proceeds just like the proof of Proposition \ref{htopy-inv}, and requires us to introduce moduli spaces of solutions to a flow that interpolates between $\xgcsql$ and $\xgcsq$. Fix $\lambda\gg 0$, and recall the function $f:(0,\infty]\to [0,1)$ from Lemma \ref{differentiability-homeomorphism}. Let $r=f(\lambda)$, and let $\alpha:\R\to [0,r]$ be a smooth non-increasing function with $\alpha(t)=r$ for $t\ll 0$, $\alpha(t)=0$ for $t\gg 0$, and let $\beta=f^{-1}\circ\alpha:\R\to[\lambda,\infty)$. Then our interpolating flow will be $t\mapsto \mc{X}^{\gc,\sigma}_{\mathfrak q^{\beta(t)}}$, which is equal to $\xagcsql$ for $t\ll0$ and is equal to $\xagcsq$ for $t\gg 0$. 

Now suppose $B^{\lambda}_-$ and $B_+$ are stationary submanifolds of $\xagcsql$ and $\xagcsq$, respectively, contained in $(W^\sigma\cap B(2R))/S^1$. We first wish to define a configuration space $\mc B_{k,\delta}(B^\lambda_-,B_+)$ of $L^2_{k,\delta}$ paths between $B^\lambda_-$ and $B_+$. Following the discussion before Lemma 3.2.11 in \cite{lin}, we can choose a smooth family of reference paths \eq{
[\gamma_0(b^\lambda_-,b_+)]:\R \to (W^\sigma\cap B(2R))/S^1} parametrized by $b^\lambda_-\in B^\lambda_-,\; b_+\in B_+$,  with the property that \eq{
    [\gamma_0(b^\lambda_-,b_+)](t)=\begin{cases}
    b^\lambda_- & t\ll 0\\
    b_+ & t\gg 0\;.\end{cases}} Given the family of paths $\gamma_0$, we then define
    \begin{multline*}
    \mc B_{k,\delta}(B^\lambda_-,B_+)=\{[\gamma]\in L^2_{k,\text{loc}}(\R,(W^\sigma\cap B(2R))/S^1)\mid [\gamma]-[\gamma_0(b^\lambda_-,b_+)]\in L^2_{k,\delta}\\
    \text{ for some } b^\lambda_-\in B^\lambda_-, b_+\in B_+\}.
    \end{multline*}
It is clear that the definition of $\mc B_{k,\delta}(B^\lambda_-,B_+)$ does not depend on the particular family of reference paths chosen. 

Now we obtain a moduli space of solutions to the interpolating flow just as in \S\ref{sec:continuation-maps}: Choose an increasing diffeomorphism $\rho:\R\to(-1,1)$, and define a lifted vector field $\tilde v$ on $(W^\sigma\cap B(2R))/S^1)\times \R$ by $\tilde v(x,s)=(\mc{X}^{\gc,\sigma}_{\mathfrak q^{\beta(s)}}(x),\rho'(s))$. Then we obtain a moduli space of trajectories of $\tilde v$:  \eq{
    M(B^\lambda_-,B_+)=\set{([\gamma],u)\in \mc B_{k,\delta}(B^\lambda_-,B_+)\times L^2_{k,\delta}(\R)\mid \frac{d}{dt}([\gamma],u)=\tilde v([\gamma],u)}.}
By arguments analogous to Theorem 3.2.12 in \cite{lin}, we have that $M(B^\lambda_-,B_+)$ is independent of $k,\delta$ when $k$ is chosen sufficiently large and $\delta$ sufficiently small.

Finally, we obtain the following proposition as in Propositions 3.6.12 and 3.6.14 in \cite{lin} (cf. also Propositions \ref{interpolating-delta-chain} and \ref{prop:codim-one-findim}):
\begin{prop}
The space $M(B^\lambda_-,B_+)$ above admits a compactification $M^+(B^\lambda_-,B_+)$ whose elements are tuples $(\bm{\gamma_0}, \gamma,\bm{\gamma_1})$, where 
\begin{enumerate}[(1)]
\item $\bm{\gamma_0}=(\gamma_0^1,\dots,\gamma_0^{n_0})$ is a broken trajectory of $\xagcsql$ beginning at $B^\lambda_-$ and ending at some intermediate $B^\lambda_+$;
\item $\bm{\gamma_1}=(\gamma_1^1,\dots, \gamma_1^{n_1})$ is a broken trajectory of $\xagcsq$ beginning at some intermediate $B_-$ and ending at $B_+$;
\item $\gamma\in M(B^\lambda_+,B_-)$.
\end{enumerate}

The compactification $M^+(B^\lambda_-,B_+)$ has the structure of an abstract $\delta$-chain whose codimension-1 faces consist of the tuples $(\bm{\gamma_0},\gamma,\bm{\gamma_1})$ above where either
\begin{enumerate}[(1)]
\item $\bm{\gamma_0}$ is unbroken $(n_0=1)$, $\bm{\gamma_1}$ is trivial $(n_1=0)$, and $\gamma_0^1$ and $\gamma$ are not boundary-obstructed;
\item $\bm{\gamma_1}$ is unbroken $(n_1=1)$, $\bm{\gamma_0}$ is trivial $(n_0=0)$, and $\gamma_1^1$ and $\gamma$ are not boundary-obstructed; or
\item Both $\bm{\gamma_0}$ and $\bm{\gamma_1}$ are unbroken $(n_0=n_1=1)$ and $\gamma_0^1$, $\gamma_1^1$ are not boundary-obstructed, but $\gamma$ is boundary-obstructed.
\end{enumerate}
\label{prop:int-mod-1}
\end{prop}

To discuss the actions of $Q$ and $V$, we will also need to consider cut-down moduli spaces. As in Section \S \ref{subsubsection:QV}, choose a generic section $\zeta$ of the canonical quaternionic line bundle over the space $(W^\sigma\cap B(2R))/\pin(2)$ and a generic section $\eta$ of the real line bundle over $(W^\sigma\cap B(2R))/\pin(2)$ corresponding to the projection $\pi:(W^\sigma\cap B(2R))/S^1\to (W^\sigma\cap B(2R))/\pin(2)$. Let $Z_\zeta$, $Z_\eta$ be the zero loci of $\zeta,\eta$, and define \eq{
    M^{\cap\zeta}(B^\lambda_-,B_+)&=\set{([\gamma],u)\in M(B^\lambda_-, B_+)\mid \pi(\gamma_0)\in Z_\zeta}\\
    M^{\cap\eta}(B^\lambda_-,B_+)&=\set{([\gamma],u)\in M(B^\lambda_-, B_+)\mid \pi(\gamma_0)\in Z_\eta}.
    }

For the cut-down moduli spaces, we have the following analogue of Proposition \ref{prop:int-mod-1}:
\begin{prop}
For a generic section $\zeta$, the cut-down moduli space $M^{\cap\zeta}(B^\lambda_-,B_+)$ admits a compactification $M^{\cap\zeta,+}(B^\lambda_-,B_+)$ whose elements are tuples $(\bm{\gamma_0}, \gamma,\bm{\gamma_1},\gamma^\ast)$, where 
\begin{enumerate}[(1)]
\item $\bm{\gamma_0}=(\gamma_0^1,\dots,\gamma_0^{n_0})$ is a broken trajectory of $\xagcsql$ beginning at $B^\lambda_-$ and ending at some intermediate $B^\lambda_+$;
\item $\bm{\gamma_1}=(\gamma_1^1,\dots, \gamma_1^{n_1})$ is a broken trajectory of $\xagcsq$ beginning at some intermediate $B_-$ and ending at $B_+$;
\item $\gamma\in M(B^\lambda_+,B_-)$.
\item $\gamma^\ast$ is one of the trajectories in the set \eq{
 \set{\gamma_0^1,\dots, \gamma_0^{n_1},\gamma,\gamma_1^1,\dots, \gamma_1^{n_1}},} and $\pi(\gamma^\ast(0))\in Z_\zeta$.
\end{enumerate}
The compactification $M^{\cap\zeta,+}(B^\lambda_-,B_+)$ has the structure of an abstract $\delta$-chain whose non-boundary-obstructed codimension-1 faces are homeomorphic to
\begin{enumerate}[(1)]
\item $M^{\cap\zeta}(B^\lambda_-, B^\lambda_+)\times_{\ev} M(B^\lambda_+,B_+)$ or $M(B^\lambda_-, B^\lambda_+)\times_{\ev} M^{\cap\zeta}(B^\lambda_+,B_+)$;
\item $M^{\cap\zeta}(B^\lambda_-,B_-)\times_{\ev} M(B_-,B_+)$ or $M(B^\lambda_-,B_-)\times_{\ev} M^{\cap\zeta}(B_-,B_+)$;
\item $M^{\cap\zeta}(B^\lambda_-,B^\lambda_+)\times_{\ev} M(B^\lambda_+,B_-)\times_{\ev} M(B_-,B_+)$ or $M(B^\lambda_-,B^\lambda_+)\times_{\ev} M^{\cap\zeta}(B^\lambda_+,B_-)\times_{\ev} M(B_-,B_+)$ or $M(B^\lambda_-,B^\lambda_+)\times_{\ev} M(B^\lambda_+,B_-)\times_{\ev} M^{\cap\zeta}(B_-,B_+)$.
\end{enumerate}
The same statements hold when the section $\zeta$ is replaced everywhere by the section $\eta$.
\label{prop:int-mod-cutdown}
\end{prop}

For an integer $M$, denote by $(\check C_\lambda)_\ast^{\le M}$, $\widecheck{\CM}_\ast^{\le M}$ the truncations of the corresponding chain complexes to degrees $|i|\le M$. Using the moduli spaces $M^+(B^\lambda_-, B_+)$, we can define a map $F: (\check C_\lambda)_\ast^{\le N-3}\to \widecheck{\CM}_\ast^{\le N-3}$ as follows. If $\sigma:\Delta\to B^\lambda_-$ is a $\delta$-chain which defines an element $[\sigma]\in (\check C_\lambda)_\ast^{\le N-3}$, then Remark \ref{rmk:negl-chains}, coupled with the fact that $B_\lambda^-$ must have dimension $0$ or $2$, shows that the index of $B_\lambda^-$ is in $[-N,N-3]$. By the same token, if $B_+$ is any stationary submanifold of $\xagcsq$, then $(\sigma\times_{\ev_-} M^+(B_-^\lambda,B_+), \ev_+)$ is nontrivial in homology if and only if $0\le \dim \sigma\times_{\ev_-} M^+(B_-^\lambda,B_+)\le 3$; but since the map $\sigma\mapsto \sigma\times_{\ev_-} M^+(B_-^\lambda,B_+)$ preserves the degree, we must have $\text{ind}(B_+)=\text{deg }\sigma - \dim \sigma\times_{\ev_-} M^+(B_-^\lambda,B_+)$, or $\deg \sigma -3 \le \ind B_+\le \deg \sigma$, which shows that $\ind B_+$ lies in $[-N,N-3]$ as well.

Let $\mathfrak C^o$, $\mathfrak C^s$, $\mathfrak C^u$ denote the set of irreducible, boundary-stable, and boundary-unstable stationary submanifolds of $\xagcsq$, and let $\mathfrak C^o_\lambda$, $\mathfrak C^s_\lambda$, $\mathfrak C^u_\lambda$ similarly denote stationary submanifolds of $\xagcsql$. For $\theta,\theta'\in \set{o,s,u}$ and a $\delta$-chain $\sigma$ in $B^\lambda_-\in\mathfrak C^\theta_\lambda$, we define \eq{
	F^\theta_{\theta'}(\sigma)=\sum_{B_+\in\mathfrak C^{\theta'}} \sigma\times_{\ev_-} M^+(B^\lambda_-,B_+).}
By the discussion above, we can assume in the above sum that $B_+$ has index in $[-N,N]$ as well, i.e. it lies in $\mc N$.
We then set \eqq{
	\check F=\begin{pmatrix}
 F^o_o& F^u_o\partial^s_u+\partial^u_oF^s_u\\
 F^o_s& F^s_s+F^u_s\partial^s_u+\partial^u_sF^s_u \end{pmatrix}:\widecheck{C}^{\lambda, \le N-3}_\ast\to \widecheck{\CM}^{\le N-3}_\ast.\label{chain-map-end}}

Now we obtain the following proposition just as in the finite-dimensional case (Proposition \ref{chain-map}):
\begin{prop}
The map $\check F$ defined in (\ref{chain-map-end}) is a chain map of degree 0.
\label{prop: chain-map-end}
\end{prop}
Now we finish the proof of Proposition \ref{prop:lambda-ultimate}.

\begin{proof}[Proof of Proposition \ref{prop:lambda-ultimate}]
Let us start by proving the statement for $\widecheck{\HM}$. We must show that the chain map $\check F$ above is invertible up to homotopy, and therefore gives an isomorphism on homology (apart from at the truncation boundary). This is entirely analogous to the proof of Proposition \ref{htopy-inv} given at the very end of Section \ref{sec:continuation-maps}: we construct a moduli space of flows on $((W\times B(2R))^\sigma)/S^1 \times \R^2$ that interpolates between the trivial flow and the composition of $\check F$ with its ``reverse'', obtained by exchanging $s$ with $-s$ in the definition of the vector field $\tilde v(x,s)$ used to construct $\check F$.

The only new point to discuss is the preservation of the $\mathcal R$-module structure. We can mimic the definition of $\check F$ using the cut-down moduli spaces, yielding \eq{
    H^\theta_{\theta'}(\sigma)&=\sum_{B_+\in\mathfrak C^{\theta'}} \sigma\times_{\ev_-} M^{\cap\zeta,+}(B^\lambda_-,B_+)\\
    \check H&=\begin{pmatrix}
 H^o_o& H^u_o\partial^s_u+\partial^u_oH^s_u\\
 H^o_s& H^s_s+H^u_s\partial^s_u+\partial^u_sH^s_u \end{pmatrix}.}

From the description of the codimension-1 stratum of $M^{\cap\zeta,+}(B^\lambda_-,B_+)$ in Proposition \ref{prop:int-mod-cutdown}, we obtain \eq{
    (\check F\circ V + V\circ \check F)+(\check H\circ \partial+ \partial\circ \check H)=0.}
Passing to homology, we conclude that $\check F$ commutes the actions of $V$ on $\widecheck{\HM}^\lambda$ and $\widecheck{\HM}$. The same argument applies to the action of $Q$.

Since $\check F$ commutes with the involution $\jmath$, this also proves the statement for the invariant $\widecheck{\HS}$.
\end{proof}
\begin{remark}
In the proof of Proposition \ref{prop:lambda-ultimate} above, we represented the $Q$ and $V$ maps as fiber products with cut-down moduli spaces. This closely mirrors the representation of the $U$ map given in Section 4.11 in \cite{kronheimerMonopolesLensSpace2007} and in Section 2.7 of \cite{lm}. In \cite{lin}, the $Q$ and $V$ maps are defined as in \cite{m3m}, using cap products with \v Cech cohomology classes on the moduli space of a cobordism. However, in this case we can replace the more general \v Cech cohomology construction with the description with cut-down moduli spaces, as the zero sets $Z_\zeta$, $Z_\eta$ used in the definition of the cut-down moduli spaces are Poincar\'e dual to the \v Cech cohomology classes in question (cf. discussion in Section \ref{subsubsection:QV}).
\label{rmk:qv-equivalence}
\end{remark}
\section{Proof of Theorem \ref{thm:main-thm}}\label{sec:main-proof}
In this section, we complete the proof of Theorem \ref{thm:main-thm}, showing that \eqq{
    \widecheck{\HS}_\ast(Y,\mathfrak s)\simeq \widetilde{H}^{\pin(2)}_\ast(\SWF(Y,\mathfrak s);\bb F) \label{eqn:restatement-main}} as modules over the ring $\mathcal R=\bb F[Q,V]/(Q^3)$.

The first point to address is the definition of the invariant $\widecheck{\HS}_\ast(Y,\mathfrak s)$. In this paper, we have defined $\widecheck{\HS}_\ast(Y,\mathfrak s)$ using configurations in global Coulomb gauge, whereas the invariant defined by Lin in \cite{lin} does not restrict to global Coulomb gauge. Nevertheless, the definition given here and the definition given in \cite{lin} are equivalent. The analysis necessary to prove this equivalence is performed in Chapter 5 of \cite{lm} in the $S^1$-equivariant setting, and carries over to the current $\pin(2)$-equivariant setting.

The object $\SWF(Y,\mathfrak s)$ on the right-hand side of (\ref{eqn:restatement-main}) is defined using the vector field $\mathcal X^{\gc}$, without any reference to a perturbation. However, it is possible to replace the vector field $\mathcal X^{\gc}$ with a perturbed vector field $\mathcal \xgcq$, producing a perturbed spectrum $\SWF_{\mathfrak q}(Y,\mathfrak s)$ which is $\pin(2)$-equivariantly stably homotopy equivalent to the unperturbed version and therefore satisfies \eqq{
    \widetilde{H}^{\pin(2)}_\ast(\SWF_{\mathfrak q}(Y,\mathfrak s);\bb F) &= \widetilde{H}^{\pin(2)}_\ast(\SWF(Y,\mathfrak s);\bb F).
    \label{eqn:pert-hom}} This construction is described in detail in Chapter 6 of \cite{lm} for the $S^1$-equivariant case, and can be exported to the current $\pin(2)$-equivariant case without change. 
    
Recall that $(\check C_\lambda,\check \partial_\lambda)$ is the Morse-Bott complex for the vector field $\xagcsql$ on the finite-dimensional manifold $(W^\lambda\cap B(2R))^\sigma/S^1$, and is defined using the isolated invariant set $\scr S^\lambda$ of points on trajectories that remain in $B(2R)$ for all time. According to Proposition \ref{prop:fin-dim-summary}, we have \eq{
    H_j(\check C_\lambda^{\text{inv}},\check \partial_\lambda)\simeq \tilde H_j^{\pin(2)}(I_{\mathfrak q}^\lambda) \quad (j\le n_\lambda-1),}
    where $I_{\mathfrak q}^\lambda$ is the Conley index for $\mathscr S^\lambda$ and $n_\lambda$ is the connectivity of $(I_{\mathfrak q}^\lambda,I_{\mathfrak q}^\lambda-(I_{\mathfrak q}^\lambda)^{S^1})$. The isomorphism respects the actions of $\mathcal R$.

However, we defined $\widecheck{\CS}^\lambda(Y,\mathfrak q)$ to be $(\check C_\lambda^{\text{inv}},\check \partial_\lambda)$ with gradings shifted down by $\dim W^{(-\lambda,0)} + 2n(Y,\mathfrak s, g)$. This grading shift agrees with the grading shift used in the definition of the spectrum $\SWF_{\mathfrak q}(Y,\mathfrak s)$. Therefore, for $\lambda=\lambda_i^{\bullet}\gg 0$, we have \eqq{
    \widecheck{\HS}^\lambda_j(Y,\mathfrak s)\simeq \widetilde{H}^{\pin(2)}_j(\SWF_{\mathfrak q}(Y,\mathfrak s);\bb F) \label{eqn:lambda-only-iso}}
for $-\dim W^{(-\lambda,0)} - 2n(Y,\mathfrak s, g) \le j \le n_\lambda-1-\dim W^{(-\lambda,0)} - 2n(Y,\mathfrak s, g)$. Setting \eq{
    M_\lambda=\min\left(\dim W^{(-\lambda,0)} + 2n(Y,\mathfrak s, g), n_\lambda-1-\dim W^{(-\lambda,0)} - 2n(Y,\mathfrak s, g)\right),} we see that (\ref{eqn:lambda-only-iso}) holds for $|j|\le M_\lambda$. Moreover, $M_\lambda\to\infty$ as $\lambda\to\infty$, as shown by equation (14.0.4) in \cite{lm}.

Combining (\ref{eqn:pert-hom}) and (\ref{eqn:lambda-only-iso}) with Proposition \ref{prop:lambda-ultimate}, we obtain the desired equivalence of $\mathcal R$-modules \eq{
    \widecheck{\HS}_\ast(Y,\mathfrak s)\simeq \widetilde{H}^{\pin(2)}_j(\SWF(Y,\mathfrak s);\bb F).}

\bibliographystyle{alpha}
\bibliography{bibliography}

\end{document}